\documentclass[10pt]{article}
\usepackage{epsf}
\usepackage{epsfig}
\usepackage{enumerate}
\usepackage{url}

\usepackage{amsmath,amsfonts,amssymb,amsthm}
\usepackage[retainorgcmds]{IEEEtrantools}

\usepackage{hyperref}
\def\pathPic{Illustrations}
\def\<{\langle}
\def\>{\rangle}

\def\sm{\setminus}
\def\nl{\newline}

\def\cB{{\cal B}}

\def\cE{{\cal E}}
\def\cF{{\cal F}}

\def\cH{{\cal H}}
\def\cI{{\cal I}}

\def\cL{{\cal L}}
\def\cM{{\cal M}}

\def\cO{{\cal O}}
\def\cP{{\cal P}}

\def\cR{{\cal R}}

\def\cT{{\cal T}}

\def\Chi{\raise .3ex
\hbox{\large $\chi$}} 
\def\lsima{\hbox{\kern -.6em\raisebox{-1ex}{$~\stackrel{\textstyle<}{\sim}~$}}\kern -.4em}
\def\lsim{\hbox{\kern -.2em\raisebox{-1ex}{$~\stackrel{\textstyle<}{\sim}~$}}\kern -.2em}
\def\({\Bigl (}
\def\){\Bigr )}
\def\({\Bigl (}
\def\){\Bigr )}

\newcommand{\be}{\begin{equation}}
\newcommand{\ee}{\end{equation}}
\newcommand{\bea}{$$ \begin{array}{lll}}
\newcommand{\eea}{\end{array} $$}
\newcommand{\bi}{\begin{itemize}}
\newcommand{\ei}{\end{itemize}}
\newcommand{\iref}[1]{(\ref{#1})}
\newtheorem{theorem}{Theorem}[section]
\newtheorem{remark}[theorem]{Remark}
\newtheorem{lemma}[theorem]{Lemma}

\newtheorem{definition}[theorem]{Definition}

\newtheorem{prop}[theorem]{Proposition}

\newtheorem*{remark*}{Remark}

\def\I{{\rm \hbox{I\kern-.2em\hbox{I}}}}
\def\P{{\rm \hbox{I\kern-.2em\hbox{P}}}}
\def\sP{{\rm \hbox{\scriptsize I\kern-.2em\hbox{\scriptsize P}}}}
\def\H{{\rm \hbox{I\kern-.2em\hbox{H}}}}
\def\R{{\rm \hbox{I\kern-.2em\hbox{R}}}}
\def\N{{\rm \hbox{I\kern-.2em\hbox{N}}}}
\def\Z{{\rm {{\rm Z}\kern-.28em{\rm Z}}}}
\def\C{{\rm \hbox{C\kern -.5em {\raise .32ex \hbox{$\scriptscriptstyle
|$}}\kern-.22em{\raise .6ex \hbox{$\scriptscriptstyle |$}}\kern .4em}}}
\def\sC{{\rm \hbox{\scriptsize \rm C\kern -.6em {\raise .4ex \hbox{$\scriptscriptstyle \scriptsize|$}}\kern-.22em{\raise .5ex \hbox{$\scriptscriptstyle \scriptsize |$}}\kern .4em}}}

\def\ra{\rightarrow}
\def\Ra{\Rightarrow}
\def\ve{\varepsilon}
\def\sq{\hfill $\diamond$\\}
\def\cO{\mathcal O}
\def\b{\mathbf}

\def\cE{\mathcal E}

\newcommand\sep{\; ; \;}% separation in sets.
\newcommand\ssep{; \,}% small separation.

\newcommand\trans{\mathrm T}

\def\cM{\mathcal M}

\def\TEq{{T_{\rm eq}}}

\def\bac{$$\left\{\begin{array}}
\def\eac{\end{array}\right.$$}
\def\beqa{\begin{eqnarray*}}
\def\eeqa{\end{eqnarray*}}

\def\gsim{\hbox{\kern -.2em\raisebox{-1ex}{$~\stackrel{\textstyle>}{\sim}~$}}\kern -.2em}

\DeclareMathOperator\Id{Id}

\DeclareMathOperator\diam{diam}

\DeclareMathOperator\disc{disc}

\DeclareMathOperator\interp{{\rm I}}

\DeclareMathOperator\Span{Span}
\def\bT{\b T}

\newcommand\seqT{(\cT_N)_{N\geq N_0}}

\def\cPi{{\text{\large$\boldsymbol\pi$}}}
\def\Lpol{\b L}
\DeclareMathOperator\M{M}
\DeclareMathOperator\GL{GL}
\DeclareMathOperator\SL{SL}

\def\TRect{{T_0}}

\def\ti{\tilde}

\def\KEq{K_*}

\newtheorem{question}{Question}
\newtheorem{programme}[question]{Programmation}

\def\bques{\begin{question}}
\def\eques{\end{question}}
\def\bprog{\begin{programme}}
\def\eprog{\end{programme}}

\usepackage{stmaryrd}

\def\II{\mathrm {I\kern-0.1exI}}

\newcommand\stext[1]{\ \text{ #1 } \ }

\DeclareMathOperator\interior{int}

\def\sR{{\rm \hbox{\scriptsize I\kern-.2em\hbox{\scriptsize R}}}}

\def\sN{{\rm \hbox{\scriptsize I\kern-.2em\hbox{\scriptsize N}}}}
\def\sC{{\rm \hbox{\scriptsize \rm C\kern -.6em {\raise .4ex \hbox{$\scriptscriptstyle \scriptsize|$}}\kern-.22em{\raise .5ex \hbox{$\scriptscriptstyle \scriptsize |$}}\kern .4em}}}

\def\CDiff{{C_\star}}
\newtheorem{thIntro}{Theorem}

\begin{document}
\title{\bf Optimally adapted meshes for finite elements of arbitrary order  
and $W^{1,p}$ norms} 
\author{Jean-Marie Mirebeau}
\maketitle
\date{}
\begin{abstract}
Given a function $f$ defined on a bounded polygonal domain $\Omega\subset \R^2$ 
and a number $N>0$, we study the
properties of the triangulation $\cT_N$ that minimizes
the distance between $f$ and its interpolation
on the associated finite element space,
over all triangulations of at most $N$ elements.
The error is studied in the $W^{1,p}$ semi-norm for $1\leq p< \infty$,
 and we consider Lagrange finite
elements of arbitrary polynomial order $m-1$.
We establish sharp asymptotic error estimates as $N\to +\infty$ when the optimal
anisotropic triangulation is used.
A similar problem has been studied in \cite{BBLS,CSX,C3, CMi2, Mi}, 
but with the error measured in the $L^p$ norm. 
The extension of this analysis to the $W^{1,p}$ norm
is required in order to match more closely the needs of numerical PDE analysis,
and it is not straightforward. In particular, the meshes
which satisfy  the optimal error estimate are characterized
by a metric describing the local aspect ratio of each triangle
and by a geometric constraint on their maximal angle, a second 
feature that does not appear for the $L^p$ error norm.
Our analysis also provides with practical strategies
for designing meshes such that the interpolation error satisfies the optimal estimate 
up to a fixed multiplicative constant. 
\vspace{-0.1cm}
\vspace{-0.5cm}
\end{abstract}

%\tableofcontents
\section*{Introduction}

In finite element approximation, a usual distinction is between {\it uniform} and
{\it adaptive} methods. In the latter, the elements defining the mesh 
may vary strongly in size and shape for a better adaptation
to the local features of the approximated function $f$. 
Such procedures are used to improve the efficiency of numerous numerical methods in scientific computing.
%This procedure allows to improve the efficiency in numerous applications of scientific computing.
This naturally raises
the objective of characterizing and constructing an {\it optimal mesh}
for a given function $f$. 
%Depending on the context, the function $f$ may be fully
%known to us, either through an explicit formula or a discrete sampling, 
%or observed through noisy measurements,
%or implicitly defined as the solution of a given
%partial differential equation. 

In this paper we consider a bounded bidimensional polygonal domain $\Omega\subset \R^2$, a fixed integer $m \geq 2$ and an exponent $1 \leq p<\infty$. 
For a given conforming triangulation $\cT$ of $\Omega$ we denote by $\interp_\cT^{m-1}$ 
the standard interpolation operator on the space of 
Lagrange finite elements of degree $m-1$ associated to $\cT$.
A general objective is to study, for any $f\in C^{m}(\overline \Omega)$, the optimization problem 
\be
\label{optW1P}
\inf_{\#(\cT)\leq N} \|\nabla (f-\interp_\cT^{m-1} f)\|_{L^p(\Omega)}
\ee
where the minimum is taken over all (possibly anisotropic) triangulations of 
cardinality $\leq N$. All the triangulations considered in this paper are assumed to be \emph{conforming}: they have no hanging nodes, which implies that the interpolant $\interp_\cT^{m-1} f$ is continuous and thus belongs to $W^{1,\infty}(\Omega)$.
The choice of the $W^{1,p}$ semi-norm appearing in the expression \iref{optW1P} is motivated by PDE analysis, e.g.\ elliptic equations in the case $p=2$. Yet our paper in mainly a contribution to approximation theory, and aims at  characterizing the approximation power of finite elements on anisotropic meshes.
We obtain sharp estimates of the \emph{asymptotical behavior} of the quantity \iref{optW1P} as $N \to \infty$, in terms of the $m$-th derivatives of $f$. 
We also describe practical strategies for constructing meshes that behave similar to the optimal one, in the sense that  they asymptotically satisfy as $N \to \infty$ this sharp error estimate up to a fixed multiplicative constant.
Obtaining similar estimates and constructions in a non asymptotic setting remains an open question.

Estimates of a similar asymptotical nature %that type
were obtained in \cite{CSX,BBLS,Loseille} in the particular case of linear
finite elements and with the error measured in the $L^p$ norm, instead of the $W^{1,p}$ semi-norm. % considered in \iref{optW1P} and throughout this paper. 
They have the form 
\be
\limsup_{N\to +\infty} \(N\min_{\#(\cT)\leq N}\|f-\interp_\cT^1 f\|_{L^p(\Omega)}\)  \leq C\left\| \sqrt{|\det(d^2f)|}\right\|_{L^{p/(p+1)}(\Omega)}, %\;\; \frac 1 \tau=\frac 1 p+1,
\label{optiaffine}
\ee
which reveals that the convergence rate is governed by the quantity $\sqrt{|\det(d^2f)|}$, which depends
nonlinearly on the Hessian $d^2f$. This is heavily tied to the fact that we allow
triangles with possibly highly anisotropic shape.
The convergence estimate \iref{optiaffine} has been extended 
to %finite elements of arbitrary order in \cite{Mi}, 
arbitrary approximation order in \cite{Mi}, where the quantity governing 
the convergence rate for finite elements of arbitrary degree $m-1$ was identified. 
This quantity depends nonlinearly on the $m$-th order derivative $d^mf$.
See also the book chapter \cite{CMi2} for an introduction to the subject
of adaptive and anisotropic piecewise polynomial approximation.
%
%The purpose of the present article is to investigate this problem when the $L^p$-norm is
%replaced by the $W^{1,p}$ semi-norm which plays a critical role in PDE analysis. This semi-norm is defined as 
%$$
%|f|_{W^{1,p}(\Omega)} := \|\nabla f\|_{L^p(\Omega)} = \left(\int_\Omega |\nabla f|^p\right)^{1/p}.
%$$
%Our second objective is to propose simple and practical ways
%of designing meshes which behave similar to the optimal one,
%in the sense that they satisfy the sharp error estimate up
%to a fixed multiplicative constant. 

\subsection*{Main results and layout}

The Taylor development of a function $f\in C^m(\overline \Omega)$, close to a point $z\in \Omega$, can be written under the form 
\be
\label{taylorMuPi}
f(z+h) = \mu_z(h)+ \pi_z(h) + o(|h|^m),
\ee
where $h\in \R^2$ is small, and where $\mu_z$ and $\pi_z$ are polynomials which respectively belong to the spaces %of polynomials
\be
\label{defPH}
\P_{m-1} := \Span\{ x^ky^l\; ; \; k+l\leq m-1\} \stext{ and } \H_m:={\rm Span}\{ x^ky^l\; ; \; k+l=m\}.
\ee
For any triangle $T$, we denote by $\interp_T^{m-1}$ the local interpolation
operator acting from $C^0(T)$ onto $\P_{m-1}$. 
For any continuous fonction $\nu \in C^0(T)$, the interpolating polynomial $\interp_T^{m-1} \nu\in \P_{m-1}$ is defined by the conditions
$$
\interp_T^{m-1} \nu(\gamma)=\nu(\gamma),
$$
for all points $\gamma\in T$ with barycentric coordinates in
the set $\{0, \frac 1 {m-1},\cdots,1\}$. 
If $T$ is a sufficiently small triangle containing the point $z$, we thus have at least heuristically on $T$
\be
\label{eqLocF}
%f-\interp_T^{m-1}f \simeq \pi_z - \interp_T^{m-1} \pi_z,
\nabla (f-\interp_T^{m-1}f) \simeq \nabla (\pi_z - \interp_T^{m-1} \pi_z),
\ee
since the Lagrange interpolation operator $\interp_T^{m-1}$ on the triangle $T$ reproduces the elements of $\P_{m-1}$. 

A key ingredient in this paper is the {\it shape function} $L_{m,p}$, which is  
defined by a {\it shape optimization problem}: 
%\begin{definition}[Shape function]
for any $\pi\in \H_m$, we define
\be
L_{m,p}(\pi):=\inf_{|T|=1} \|\nabla (\pi - \interp_T^{m-1} \pi)\|_{L^p(T)},
\label{shapeFunctionL}
\ee
where the infimum is taken over all triangles of area $|T|=1$. 
%\end{definition}
The solution to this optimization problem thus describes the shape of the triangles 
of area $1$ which are best adapted
to the polynomial $\pi$ in the sense of minimizing the interpolation error
measured in $W^{1,p}$. 
If $\cT$ is a triangulation of a domain $\Omega$, then $\interp_\cT^{m-1}$ 
refers to the interpolation operator which coincides with $\interp_T^{m-1}$ on each %triangle 
$T\in \cT$. In view of \iref{eqLocF}, the optimization problem appearing in \iref{shapeFunctionL} can be regarded as a ``local'' version of the ``global'' problem \iref{optW1P} of interest.

The function $L_{m,p}$ is the natural generalisation of the function $K_{m,p}$ introduced in \cite{Mi} for the study of optimal anisotropic triangulations in the sense of the $L^p$ interpolation error: for all $\pi \in \H_m$
$$%\be
K_{m,p}(\pi):=\inf_{|T|=1} \|\pi - \interp_T^{m-1} \pi\|_{L^p(\Omega)}.
%\label{shapeFunctionK}
$$%\ee

Throughout this paper we denote by $\tau\in (0, \infty)$ the exponent defined by 
\be
\label{defTau}
\frac 1 \tau := \frac{m-1} 2 + \frac 1 p.
\ee
Consider two triangles $T$, $T'$ and a polynomial $\pi \in \H_m$. If $T$ is mapped onto $T'$ by a transformation of the form $z\mapsto \alpha z + \beta$, where $\alpha \in \R\sm \{0\}$ and $\beta\in \R^2$ (in other words the composition of a translation, an homothety, and a central symmetry if $\alpha<0$), then recalling that $\pi$ is $m$-homogeneous one %obtains 
easily checks that 
\be
\label{transInv}
|T|^{-\frac 1 \tau} \|\nabla (\pi - \interp_T^{m-1} \pi)\|_{L^p(T)} =  |T'|^{-\frac 1 \tau} \|\nabla (\pi - \interp_{T'}^{m-1} \pi)\|_{L^p(T')}.
\ee
Therefore $\|\nabla (\pi-\interp_T^{m-1} \pi)\|_{L^p(T)} \geq |T|^\frac 1 \tau L_{m,p}(\pi)$ for any triangle $T$ and any $\pi \in \H_m$.
Our asymptotic error estimate for the optimal triangulation is given by
the following theorem. 
\begin{thIntro}
\label{mainTheorem}
For any bounded polygonal domain $\Omega\subset \R^2$, any function 
$f\in C^m(\overline \Omega)$, and any $1\leq p <\infty$, there exists a sequence 
of triangulations $\seqT$ of $\Omega$, with $\#(\cT_N)\leq N$, such that
\be
\label{upperEstim}
\limsup_{N\ra \infty} N^{\frac {m-1} 2} \|\nabla (f-\interp_{\cT_N}^{m-1} f)\|_{L^p(\Omega)} \leq 
\left\|L_{m,p}\left(\frac{d^m f}{m!}\right)\right\|_{L^\tau(\Omega)}. %,\text{ where } \frac 1 \tau := \frac {m-1} 2 +\frac 1 p 
\ee
\end{thIntro}
This theorem is the consequence of a sharper result, Theorem \ref{optiTheorem}, which is given below.
In the above estimate we slightly abuse notations by defining for each $z\in \Omega$
\be
\label{defIdent}
L_{m,p}\left(\frac{d^m f(z)}{m!}\right) := L_{m,p}(\pi_z),
\ee
where the polynomial $\pi_z\in \H_m$ is defined by \iref{taylorMuPi}. 
In other words we identify the collection $d^m f(z)$ of $m$-th derivatives of $f$ at a given point $z\in \Omega$ to the corresponding term in the Taylor development of $f$ close to $z$. Explicitly the right hand side of \iref{upperEstim} stands for 
$$
\left\|L_{m,p}\left(\frac{d^m f}{m!}\right)\right\|_{L^\tau(\Omega)} := \left(\int_\Omega L_{m,p}(\pi_z)^\tau dz\right)^\frac 1 \tau.
$$

The integer $N_0$ appearing in Theorem \ref{mainTheorem} is independent of $f$
and refers to the 
minimal cardinality of a conforming triangulation of $\Omega$.  
An important feature of the estimate \iref{upperEstim} is the ``$\limsup$''. 
Recall that the upper limit of a sequence $(u_N)_{N\geq N_0}$ is defined by 
$$
\limsup_{N\to \infty} u_N := \lim_{N\to \infty} \sup_{n\geq N} u_n,
$$
and is in general stricly smaller than the supremum $\sup_{N\geq N_0} u_N$. It is still an open question to find an appropriate upper estimate of $\sup_{N\geq N_0} N^{\frac{m-1} 2}  \|\nabla (f-\interp_{\cT_N}^{m-1} f)\|_{L^p(\Omega)}$ when optimally adapted anisotropic triangulations are used.\\

We show in \S \ref{secPractical} that a triangulation satisfies the optimal estimate of Theorem \ref{mainTheorem}, up to a fixed multiplicative constant, if it obeys the following four general principles:

\begin{tabular}{ll}
(i) &The interpolation error should be evenly distributed on all triangles.\\
(ii) &The triangles should adopt locally a specific aspect ratio, dictated by the local value of $d^m f$.\\ 
(iii) &The largest angle of the triangles should be bounded away from $\cPi = 3.14159\ldots$\\
(iv) &The triangulation $\cT$ should be sufficiently refined in order to adapt to the local features of $f$.\\
\end{tabular}

The third point (iii) is the main new ingredient of this paper compared to \cite{Mi}, and is necessary for obtaining optimal $W^{1,p}$ error estimates (but not for $L^p$ error estimates). Roughly speaking, two triangles having the same optimized aspect ratio imposed
by (ii) may greatly differ in term of their largest angle, and the most acute triangle should
be preferred when error is measured in $W^{1,p}$ rather than $L^p$.
The influence of large angles in mesh adaptation has already been studied in \cite{Ba,Ja,Shew2,C4}. 
The heuristic guideline is that large angles should be avoided in general, since they lead to oscillations of the gradient of the interpolant. On the contrary, extremely thin triangles and very small angles can be necessary for optimal mesh adaptation.  

%We describe in \S\ref{secPractical} the practical construction of a sequence of triangulations which asymptotically satisfies the optimal estimate \iref{upperEstim} up to a fixed multiplicative constant.
The shape function $L_{m,p}$ plays an important role in our results, and we therefore devote \S \ref{secShape} to its study which is based on algebraic techniques. 
We obtain explicit minimizers, up to a fixed multiplicative constant, of the optimization problems which correspond to piecewise linear and piecewise quadratic finite element approximation. We also introduce, for arbitrary $m\geq 2$, explicit functions $\pi\in \H_m\mapsto \Lpol_m(\pi)$ which are defined as the root of a polynomial in the coefficients of $\pi$, and are \emph{uniformly equivalent} to the shape function $L_{m,p}$, leading therefore to asymptotic error estimates similar to \iref{upperEstim} up to multiplicative constants.\\

In order to illustrate the sharpness of Theorem \ref{mainTheorem}, we introduce
a slight restriction on sequences of triangulations, following 
an idea in \cite{BBLS}: a sequence $\seqT$ of triangulations is said to be \emph{admissible} if
$\#(\cT_N) \leq  N$ and $\sup_{N \geq N_0} (N^\frac 1 2 \sup_{T\in \cT_N} \diam(T)) < \infty$, %where $\diam(E)$ denotes the diameter of a set $E\subset \R^2$. I
in other words if
\be
\label{admissibilityCond}
\sup_{T\in \cT_N} \diam(T) \leq C_AN^{-\frac 1 2}
\ee
for some constant $C_A>0$ independent of $N$. 
Here and below we denote the diameter of a set $E\subset \R^2$ by $\diam(E) := \sup \{|x-y| \ssep x,y\in E\}$. The following theorem shows that the estimate
\iref{upperEstim} cannot be improved when we restrict our attention to admissible sequences of triangulations.
It also shows that this class is reasonably large in the sense that 
\iref{upperEstim} is ensured to hold up to small perturbation.

\begin{thIntro}
\label{optiTheorem}
Let $\Omega\subset \R^2$ be a bounded polygonal domain, let $f\in C^m(\overline \Omega)$ and let $1\leq p <\infty$.
For any \emph{admissible} sequence $\seqT$ of triangulations of $\Omega$, one has
\be
\label{lowerEstim}
\liminf_{N\ra \infty} N^{\frac {m-1} 2} \|\nabla (f-\interp_{\cT_N}^{m-1} f)\|_{L^p(\Omega)} \geq \left\|L_{m,p}\left(\frac{d^m f}{m!}\right)\right\|_{L^\tau(\Omega)}.
\ee
Furthermore, for all $\ve>0$ there exists an \emph{admissible} sequence of triangulations $(\cT_N^\ve)_{N\geq N_0}$ such that
\be
\label{upperEstimEps}
\limsup_{N\ra \infty} N^{\frac {m-1} 2}  \|\nabla(f-\interp_{\cT_N^\ve}^{m-1} f)\|_{L^p(\Omega)} \leq \left\|L_{m,p}\left(\frac{d^m f}{m!}\right)\right\|_{L^\tau(\Omega)}+\ve.
\ee
\end{thIntro}
Note that the sequences of triangulations $(\cT_N^\ve)_{N\geq N_0}$ satisfy the admissibility condition \iref{admissibilityCond} with a constant $C_A(\ve)$ which may grow to $+\infty$ as $\ve\to 0$. Theorem \ref{mainTheorem} can be inferred from the estimate \iref{upperEstimEps} proceeding as follows: for each $N \geq N_0$ and each $\ve>0$ we define a real $\delta(N,\ve)$ by the equality
$$
N^{\frac {m-1} 2}  \|\nabla(f-\interp_{\cT_N^\ve}^{m-1} f)\|_{L^p(\Omega)} = \left\|L_{m,p}\left(\frac{d^m f}{m!}\right)\right\|_{L^\tau(\Omega)}+ \delta(N, \ve).
$$
We next observe that for any fixed $\ve>0$ one has $\limsup_{N\to \infty} \delta(N,\ve) \leq \ve$. We may therefore construct, using a diagonal extraction procedure, a sequence $(\ve_N)_{N \geq N_0}$ such that $\limsup_{N \to \infty}\delta(N, \ve_N) \leq 0$. The sequence of triangulations $(\cT_N^{\ve_N})_{N \geq N_0}$ then clearly satisfies \iref{upperEstim}, which establishes Theorem \ref{mainTheorem}.

The proof of Theorem \ref{optiTheorem} is given in \S \ref{secProofs}.
The proof of the upper estimate \iref{upperEstimEps} involves the construction
of an optimal mesh based on a patching strategy adapted from the one encountered in \cite{BBLS}. 
However, inspection of the proof reveals that
this construction only becomes effective as
the number of triangles $N$ becomes very large. The construction described in \S \ref{secPractical} should therefore be preferred in practical applications.\\

%\begin{remark} 
%It can easily be shown that if $(\cT_N)_{N\geq N_0}$ is an admissible sequence of triangulations and
%$f\in\cC^m(\Omega)$, then $\|f-\interp_{\cT_N}^m f\|_{L^p(\Omega)}$
%decays with the rate $N^{-m/2}$ which is faster than the 
%decay rate $N^{-\frac{m-1} 2}$ obtained for the $W^{1,p}$ error. The convergence estimate \iref{upperEstimEps}
%is therefore also valid in the $W^{1,p}$ norm 
%$$
%\|f\|_{W^{1,p}(\Omega)} := \(\|f\|_{L^p(\Omega)}^p+|f|_{W^{1,p}(\Omega)}^p\)^{1/p}.
%$$ 
%\end{remark}

%We discuss in \S \ref{secPractical} the practical construction of a triangulation which satisfies the optimal estimate \iref{upperEstim} up to a fixed multiplicative constant. It suffices to obey four general principles:

\subsection*{Notations}
%For any pair of vectors $u,v\in\R^2$,
We denote by $\<u,v\>$ the inner product of two vectors $u,v\in \R^2$, and by
$
|u|:=\sqrt{\<u,u\>}
$
the euclidean norm of $u$. When $g\in L^p(\Omega, \R^2)$ is a vector valued function,
we denote by $\|g\|_{L^p(\Omega)}$ the $L^p$ norm of $x\mapsto |g(x)|$ on $\Omega$, for instance in \iref{optW1P} and \iref{shapeFunctionL}.

We denote by $\M_2$ the set of all $2\times 2$ real matrices, equipped with the spectral norm
$
\|A\|:=\max_{|u|\leq 1}|A u|.
$
We denote by $\GL_2\subset \M_2$ the group of invertible matrices, by $\SL_2\subset \GL_2$ the special group of matrices of determinant $1$, and by $\cO_2\subset \GL_2$ the group of orthogonal matrices.
We denote by $S_2\subset \M_2$ the linear space of symmetric matrices, by $S_2^\oplus \subset S_2$ the subset of non-negative symmetric matrices, and by $S_2^+\subset S_2^\oplus$ the subset of positive definite symmetric matrices.

For any two symmetric matrices $S,S'\in S_2$, we write $S\leq S'$ if and only if $S'-S\in S_2^\oplus$.
For any $S\in S_2^\oplus$ (resp. $S_2^+$) and any $\alpha >0$ (resp. $\alpha \in \R$) we denote by $S^\alpha$ 
the symmetric matrix obtained by elevating the eigenvalues to the power $\alpha$ in a diagonalization of $S$.
 
%We denote by $S_d\subset \M_d(\R)$ the subset of symmetric matrices, and by $S_d^+\subset S_d$ the subset of non-negative symmetric matrices.
%$$
%S_d := \{S\in \M_d(\R)\sep S^\trans = S\} \text{ and } S_d^+ := \{S\in S_d \sep z^\trans S z\geq 0 \text{ for all }  z\in \R^d\}.
%$$ %\{S\in S_d \sep \forall z\in \R^d, \ z^\trans S z\geq 0\}.
%For any two symmetric matrices $S,S'\in S_d$, we write $S\leq S'$ if and only if $S'-S\in S_d^+$.
%
%We equip the spaces $\P_m$ and $\H_m$ with the norm
%\be
%\label{normhm}
%\|\pi\|:=\max_{|u|\leq 1}|\pi(u)|.
%\ee
The greek letter $\pi$ always refers to an homogeneous polynomial $\pi\in \H_m$, 
while the bold notation $\cPi$ refers to the mathematical constant $\cPi= 3.14159\ldots$

%Recall that if $g$ is a $C^m$ function, we identify $d^mg(x)$ to a polynomial in $\H_m$, by \iref{dmfHm}. We
%then denote
%\be
%\|d^m g\|_{L^\infty(E)}:=\max_{z\in E} \|d^m g(z)\|
%\label{normdmg}
%\ee
%with $\|\cdot\|$ the previously defined norm on $\H_m$.

%\section{The shape function $L_{m,p}$ and local error estimates}

\section{Adaptive mesh generation} %Local error estimates}%Optimal meshes}%Practical mesh generation}
\label{secPractical}

This section describes some properties that are needed to ensure that a mesh $\cT$ satisfies the optimal error estimate introduced in Theorem \ref{mainTheorem} up to a fixed multiplicative constant. 
We first introduce a description of the triangle based on some parameters adapted to our purposes.
%These parameters should follow a
We then prescribe the behavior of these parameters in a triangulation tailored to the approximation of a given function $f$, and we discuss the practical challenges encountered in the construction of such a mesh.
%We then show how to these parameters should discuss the practical challenges encountered in the construction of such triangulations.

%We discuss in this section the practical issues encountered when attempting to design adaptive anisotropic finite element meshes which satisfy the optimal estimate, given in Theorem \ref{mainTheorem}, up to a fixed multiplicative constant.

%We focus in this section on practical issues when designing adaptive anisotropic finite element meshes. 

\subsection{Description of a triangle}

A triangle $T\subset \R^2$ is determined by the collection of its three vertices $v_1$, $v_2$ and $v_3$. In the context of adaptive mesh generation, we rather adopt the following parameters. The position of $T$ is determined by its barycenter
$$
z_T := (v_1+v_2+v_3)/3.
$$ 
The area, the aspect ratio and the orientation of $T$ are encoded in a symmetric positive definite matrix $\cH_T\in S_2^+$ defined by the equality
%For each triangle $T$, of vertices $v_1$, $v_2$ and $v_3$, we denote by $z_T := (v_1+v_2+v_3)/3$ its barycenter, and we define a matrix $\cH_T\in S_2^+$ by the equality 
$$
\cH_T^{-1} := \frac 2 3 \sum_{1\leq i \leq 3} (v_i-z_T) (v_i-z_T)^\trans.
$$
Last we shall introduce below a real $S(T) \geq 1$ which is tied to the largest angle of $T$.

\begin{figure}
\centering
\begin{tabular}{cc}
{\raise 18mm \hbox{
\begin{tabular}{c}
%\vspace{15mm}
\includegraphics[width=5cm,height=1.6cm]{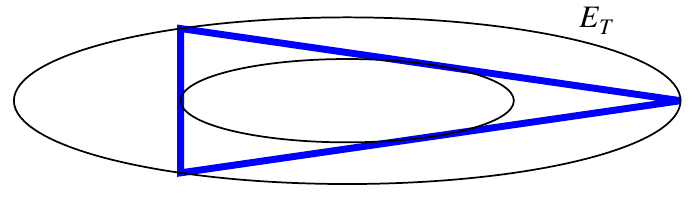}\\
%\vspace{15mm}
\includegraphics[width=5cm,height=1.6cm]{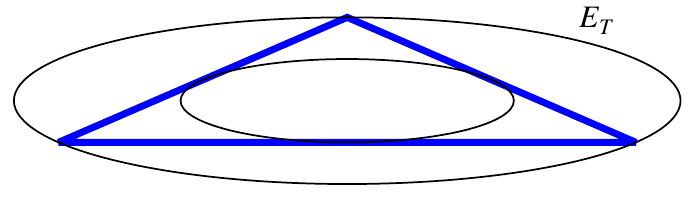}
\end{tabular}
}}
&
\includegraphics[width=3.2cm,height=4cm]{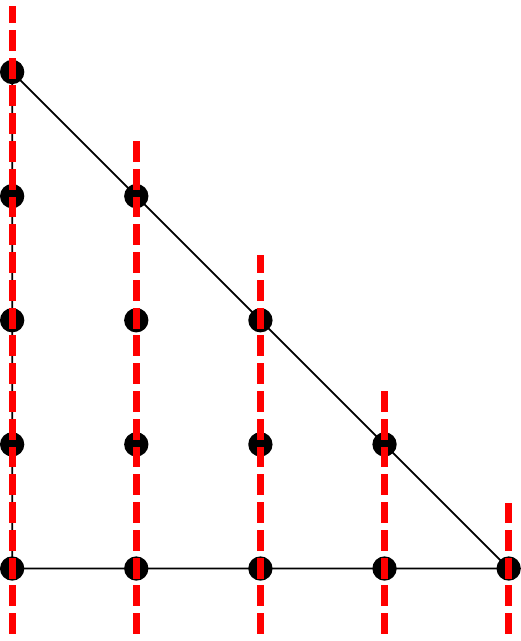} %{Illustrations/Lagrange4Lines.eps}
%&
%{\raise 5mm \hbox{
%\includegraphics[width=4cm,height=4cm]{PaperW1P/Lagrange4.pdf} %{Illustrations/Lagrange4.eps}
%}}
%%\hspace{1cm}
%&
%\includegraphics[width=4cm,height=5cm]{PaperW1P/Lagrange4Lines.pdf} %{Illustrations/Lagrange4Lines.eps}
\end{tabular}
	\caption{Two triangles and the associated ellipses associated by \iref{defET} (left). The interpolation points on the triangle $T_0$ are aligned vertically (right).}
	\label{fig1W1P}
\end{figure}

If a triangle $T'$ is mapped onto $T$ by the change of coordinates $z\mapsto Az+z_0$, where $A\in \GL_2$ and $z_0\in \R^2$, then 
%If $A$ is an invertible $2\times 2$ matrix and if $T'$ is mapped onto $T$ by the linear map $z\mapsto Az$, then 
one easily checks that 
\be
\label{eqInvHT}
\cH_{T'} = A^\trans \cH_T A.
\ee
From this point onwards we denote by $\TEq$ the triangle of vertices $(\cos(2 k  \cPi/3), \sin (2 k \cPi/3))_{0\leq k \leq 2}$, which satisfies $\cH_\TEq = \Id$. Combining this observation with \iref{eqInvHT}, Proposition 5.1.3 in \cite{thesisJM} establishes that for any triangle $T$
\be
\label{eqAreaHT}
|T|\sqrt{\det \cH_T} = |\TEq|, 
\ee
and that there exists a rotation $U\in \cO_2$ (depending on $T$) such that the change of coordinates
\be
\label{eqMapTTeq}
z\mapsto U \cH_T^\frac 1 2 (z-z_T)
\ee
maps $T$ onto $\TEq$. 
%where $\TEq$ is the triangle of vertices $(\cos(2 k  \cPi/3), \sin (2 k \cPi/3))_{0\leq k \leq 2}$.
Furthermore, as illustrated on Figure \ref{fig1W1P} (left), we have the inclusions
\be
\label{defET}
%\{z\in \R^2 \ssep (z-z_T)^\trans \cH_T (z-z_T)\leq 1/4 \} \subset T \subset \cE_T := \{z\in \R^2 \ssep (z-z_T)^\trans \cH_T (z-z_T)\leq 1\}.
\{z_T+u \ssep u^\trans \cH_T u\leq 1/4 \} \subset T \subset \cE_T := \{z_T+u \ssep u^\trans \cH_T u\leq 1\}.
\ee
The inclusion of two triangles $T,T'$ therefore implies an inclusion of ellipses, hence an inequality on the associated symmetric matrices: %(*details*)and Note that for any two triangles $T,T'$,
\be
\label{ineqTH}
T'\subset T \ \Ra \  4\cH_{T'} \geq \cH_T.%\{z_{T'}+u \ssep u^\trans \cH_{T'} u\leq 1/4 \} \subset \cE_T \ \Ra \
\ee
%The inclusion of two triangles imply the inclu
%Furthermore the ellipse of minimal volume containing $T$ is $\cE_T := \{z \ssep (z-z_T)^\trans \cH_T (z-z_T)\leq 1\}$. The matrix $\cH_T$ thus encodes the area, the aspect ratio and the orientation of $T$.
We denote by $\rho(T)\in [1, \infty[$ the \emph{measure of degeneracy} of a triangle $T$, which is defined as follows: 
\be
\label{defrho}
\rho(T) := \sqrt{\|\cH_T\| \|\cH_T^{-1}\|}.
\ee

\begin{lemma}
\label{lemmaDiamRho}
For each triangle $T$, one has %$\diam(T)^2 \leq (4 /\sqrt 3) |T| \rho(T)$.
$$
\rho(T)\leq \frac{\diam(T)^2}{|T|/|\TEq|} \leq 4\rho(T).
$$
\end{lemma}

\begin{proof}
It follows from \iref{eqAreaHT} that $|T|/|\TEq| = \|\cH_T\|^{-\frac 1 2} \|\cH_T^{-\frac 1 2}\|$, and from \iref{defET} that $\|\cH_T^{-\frac 1 2}\| = \diam(\cE_T)/2\leq \diam(T)\leq \diam(\cE_T)$.
%One easily checks that this equality holds for the equilateral triangle $\TEq$. For a general triangle $T$ we obtain using \iref{eqMapTTeq} that $\diam(\TEq) \geq \|\cH_T^{-1}\|^{-\frac 1 2}\diam(T)$, and from \iref{eqAreaHT} that 
%%$|T| \|\cH_T\|^\frac 1 2 \|\cH_T^{-1}\|^{-\frac 1 2} = |\TEq|$.
%$ |\TEq|=|T| \sqrt{\|\cH_T\| \|\cH_T^{-1}\|^{-1}}$. 
Combining this with \iref{defrho} we obtain the announced result.
\end{proof} %the above definition of $\rho(T)$

As illustrated on Figure \ref{fig1W1P}, the fact that a triangle $T$ is acute, or not, is not reflected on the ellipsoid $\cE_T$ or the matrix $\cH_T$. Since acute triangles play a priviledged role in finite element approximation %we denote their collection by $\mA$, and 
we introduce the \emph{measure of sliverness} $S(T)$ of a triangle $T$, which is defined as follows
\be
\label{defS}
S(T) := \inf\{ \|\psi\| \, \|\psi^{-1}\| \ssep \psi\in \GL_2 %\text{ s.t. } 
\text{ s.t. the image of } T \text{ by } \psi \text{ is acute}\}.
\ee
(Where we refer to the image of $T$ by the linear change of coordinates $z \mapsto \psi z$ on $\R^2$.)
The quantity $S(T)$ can be regarded as the distance from $T$ to the collection of acute simplices. It immediately follows from \iref{eqMapTTeq} and \iref{defrho} that 
\be
\label{ineqRhoS}
1 \leq S(T) \leq \rho(T).
\ee
Then next expression gives an explicit expression of $S(T)$ in terms of the largest angle of $T$, which shows that it is equivalent up to a multiplicative constant to the quantities previously introduced in \cite{ApKu,Ja} and referred to  as $\sigma_{\text{min}}(T)$ and $1/\cos(\theta')$ respectively.

\begin{prop} For any triangle $T$ with largest interior angle $\theta$, one has 
\label{propSTan}
$S(T) = \max\{1, \tan \frac \theta 2\}.$
\end{prop}

\begin{proof}
The result of this proposition is trivial if the triangle $T$ is acute, we therefore assume that $T$ is obtuse. We assume without loss of generality that the vertices of $T$ are $0$, $\alpha u$ and $\beta v$, where $\alpha, \beta>0$, $u,v\in \R^2$, $|u|=|v|=1$ and $\<u,v\> = \cos \theta$. Note that $|u-v|= 2 \sin (\theta/2)$ and $|u+v|= 2 \cos (\theta/2)$. 
Let $\psi\in \GL_2$ be such that the image of $T$ by $\psi$ is acute. We thus have $\<\psi(u),\psi(v)\>\geq 0$ and therefore $|\psi(u)-\psi(v)| \leq |\psi(u)+\psi(v)|$. It follows that 
$$
\|\psi\|\; \|\psi^{-1}\|\geq\frac {|u-v|} {|u+v|} \times \frac {|\psi(u)+\psi(v)|} {|\psi(u)-\psi(v)|} \geq \frac {2 \sin (\theta/2)}{2\cos(\theta/2)}  = \tan \frac \theta 2.
$$ 
Therefore $S(T)\geq \tan \frac \theta 2$.
On the other hand, let $\psi$ be defined by $\psi(u)=(0,1)$ and $\psi(v) = (1,0)$. Obviously the image of $T$ by $\psi$ has one of its angles equal to $\cPi/2$, and is therefore acute. One easily checks that $\|\psi\| \|\psi^{-1}\| = \tan (\theta/2)$ and therefore
$S(T)\leq \tan \frac \theta 2$, which concludes the proof of this proposition.
\end{proof}

Most error estimates available in the literature, such as in \cite{Ba,Ja}, are designed to control the gradient interpolation error of a function on a triangle, by second or higher derivatives of the approximated function.
Our purposes require a slight variant of these estimates, given in the next lemma, in which the gradient interpolation error is controlled by the gradient itself of the approximated function.
%, and are therefore not completely adapted to our purposes. The next estimate controls this error by the gradient itself of the function.

\begin{lemma}
\label{lemmaS}
There exists a constant $C=C(m)$ such that the following holds. For any triangle $T$ and any $f\in W^{1,\infty}(T)$, one has
$$
\|\nabla(f-\interp_T^{m-1} f)\|_{L^\infty(T)} \leq C S(T) \|\nabla f\|_{L^\infty(T)}.
$$
\end{lemma}

\begin{proof}
Let $T_0$ be the triangle of vertices $(0,0)$, $(1,0)$ and $(0,1)$, and let $g\in W^{1,\infty}(T_0)$. We define $\ti g(x,y):= g(x,0)$ and $h(x,y):=g(x,y)-g(x,0)$. 
Since $\ti g$ does not depend on $y$ and since the Lagrange interpolation points on $T_0$ are aligned vertically, as illustrated on Figure \ref{fig1W1P}, the Lagrange interpolant $\interp_{T_0}^{m-1} \ti g$ does not depend on $y$ either. Futhermore, for all $(x,y)\in T_0$, we have 
$|h(x,y)| = |\int_{s = 0}^y \frac {\partial g}{\partial y}(x,s) ds | \leq \|\frac {\partial g}{\partial y}\|_{L^\infty(T_0)}$.
Hence
$$
\left\|\frac {\partial \interp_\TRect^{m-1} g} {\partial y}   \right\|_{L^\infty(\TRect)}
= \left\| \frac {\partial \interp_\TRect^{m-1} h}{\partial y}  \right\|_{L^\infty(\TRect)}\\
\leq  C_0 \|\interp_\TRect^{m-1} h\|_{L^\infty(\TRect)}\\
\leq  C_0 C_1 \|h\|_{L^\infty(\TRect)} \\
\leq  C_0 C_1  \left\|\frac {\partial g}{\partial y}\right\|_{L^\infty(\TRect)},
$$
%\begin{eqnarray*}
%\left\|\frac {\partial \interp_\TRect^{m-1} g} {\partial y}   \right\|_{L^\infty(\TRect)}
%&=& \left\| \frac {\partial \interp_\TRect^{m-1} h}{\partial y}  \right\|_{L^\infty(\TRect)}\\
%&\leq & C_0 \|\interp_\TRect^{m-1} h\|_{L^\infty(\TRect)}\\
%&\leq & C_0 C_1 \|h\|_{L^\infty(\TRect)} \\
%&\leq & C_0 C_1  \left\|\frac {\partial g}{\partial y}\right\|_{L^\infty(\TRect)},
%\end{eqnarray*}
where the constants $C_0$ and $C_1$ are the $L^\infty(T_0)$ norms of the operators $g\in \P_{m-1}\mapsto \frac{\partial g}{\partial y}\in \P_{m-2}$, and $g\in C^0(T_0)\mapsto \interp_{T_0}^{m-1} g\in \P_{m-1}$ respectively.

Let $e$ be an edge vector of the triangle $T$. There exists an affine change of coordinates $\Psi$ on $\R^2$, with linear part $\psi\in \GL_2$, such that $\Psi(T) = T_0$ and $\psi e = e_0$, where $e_0=(0,1)$ is the vertical edge vector of $T_0$.
Noticing that 
$$
\<e,\nabla \interp_T^{m-1} (g\circ \Psi)\> = \<e,\nabla ((\interp_\TRect^{m-1} g)\circ \Psi)\> = \<e_0, (\nabla \interp_\TRect^{m-1} g) \circ \Psi\> = \frac{\partial  \interp_\TRect^{m-1} g}{\partial y} \circ \Psi, 
$$
we obtain
\be
\label{edgeGrad}
\|\<e, \nabla \interp_T^{m-1} (g\circ \Psi)\>\|_{L^\infty(T)} = \left\|\frac{\partial \interp_{T_0}^{m-1} g}{\partial y}  \right\|_{L^\infty(\TRect)} \\
\leq 
C_0C_1 \left\|\frac{\partial g}{\partial y} \right\|_{L^\infty(\TRect)}\\
=C_0C_1 \|\<e, \nabla (g\circ \Psi)\>\|_{L^\infty(T)}.
\ee
%\be
%\label{edgeGrad}
%\begin{array}{rcl}
%\|\<e, \nabla \interp_T^{m-1} (g\circ \Psi)\>\|_{L^\infty(T)} &=& \left\|\frac{\partial \interp_{T_0}^{m-1} g}{\partial y}  \right\|_{L^\infty(\TRect)} \\
%&\leq& 
%C_0C_1 \left\|\frac{\partial g}{\partial y} \right\|_{L^\infty(\TRect)}\\
%&=&C_0C_1 \|\<e, \nabla (g\circ \Psi)\>\|_{L^\infty(T)}.
%\end{array}
%\ee
Applying this inequality to $g=f\circ \Psi^{-1}$ we obtain that 
\be
\label{ineqSfe}
\|\<e, \nabla \interp_T^{m-1} f\>\|_{L^\infty(T)} \leq C_0C_1 \|\<e, \nabla f\>\|_{L^\infty(T)},
\ee
for any edge vector $e\in \{a,b,c\}$ of $T$. 
Defining the norm 
%$
%|v|_T := \max\{|a|^{-1} |\<a,v\>|,|b|^{-1} |\<b,v\>|,|c|^{-1} |\<c,v\>|\},
%$
$
|v|_T := \max\{|\<v,a\>|/|a|,\ |\<v,b\>|/|b|, \ |\<v,c\>|/|c|\},
$
%We next define a norm on $\R^2$ as follows
%$$
%|v|_T := |a|^{-1} |\<a,v\>|+|b|^{-1} |\<b,v\>|+|c|^{-1} |\<c,v\>|.
%$$
we obtain %from inequality \iref{ineqSfe} that 
\be
\label{gradNormT}
\| \ |\nabla \interp_T^{m-1} f|_T\ \|_{L^\infty(T)} \leq C_0C_1 \| |\nabla f|_T \|_{L^\infty(T)}.
\ee
For any $v\in \R^2$ one has 
%We next observe that if $\theta$ denotes the maximal angle of $T$, 
$
\cos(\theta/2) |v| \leq |v|_T\leq  |v|,
$
where $\theta$ denotes the maximal angle of $T$. Indeed
%where $|\cdot|$ is the euclidean norm: 
the upper inequality is trivial
and the lower one is follows from the fact that at least one of the edge vectors makes
an angle less than $\theta/2$ with $v$. Combining this with \iref{gradNormT}, we obtain
$$
\| \nabla \interp_T^{m-1} f\ \|_{L^\infty(T)} \leq \frac {C_0C_1}{\cos(\theta/2)}
 \| \nabla f \|_{L^\infty(T)}.
$$
Since $\theta \geq \pi/3$ we have $\frac 1 {\cos(\theta/2)}\leq 2 \tan(\theta/2) \leq 2S(T)$
according to Proposition \ref{propSTan}, which concludes the proof with $C=2C_0C_1+1$.
\end{proof}

\begin{remark*}
The definition \iref{defS} of the measure of sliverness $S(T)$, Propositions \ref{propSTan} and Lemma \ref{lemmaS} have analogs for simplices of arbitrary dimension, see for instance the definition (3.87), Proposition 3.6.1 and Lemma 3.6.3 in \cite{thesisJM}.
\end{remark*}

\subsection{Construction of a triangulation adapted to a function} %Adaptive mesh generation} %Local error estimates and  Optimal shape of a triangle}
\label{subsecConst}
\label{secPractMesh}
%Our first objective is to identify the optimal shape of the triangles. For that purpose we reformulate the optimization problem 

In order to identify the optimal shape of the triangles, our first step is to reformulate the optimization problem appearing in the definition \iref{shapeFunctionL} of the shape function $L_{m,p}$.
%We examine the optimization problem \iref{shapeFunctionL} which defines the shape function $L_{m,p}(\pi)$ and w
%One easily obtains, % Our first observation is that, 
%Since all polynomials $\pi\in \H_m$ are homogeneous, the problem \iref{shapeFunctionL} defining $L_{m,p}(\pi)$, the minimization of the interpolation error among all triangles of area $1$, can be cast into optimizing %reformulated into an optimization of 
%the area among all triangles of interpolation error $\leq 1$. 
One easily checks using the invariance property \iref{transInv} that for any $\pi \in \H_m$
%Note that due to the homogeneity of $\pi$, %Or equivalently 
\be
\label{optArea}
L_{m,p}(\pi) = \inf \{ |T|^{-\frac {m-1} 2} \ssep T \text{ s.t. } |T|^{-\frac 1 p} \| \nabla (\pi-\interp_T^{m-1} \pi)\|_{L^p(T)} \leq 1\},
\ee
and the minimizers of this optimization problem and the original one are homothetic.

The next lemma shows that the exponent $p$ can be disregarded in the above expression, when one is only interested in minimizing it up to a fixed multiplicative constant. %since we are only interested in minimizing \iref{optArea} up to a fixed multiplicative constant.
\begin{lemma}
\label{lemmaEquiv}
There exists a constant $c=c(m)>0$ such that for any triangle $T$ and any $\pi \in \H_m$ %and any $1 \leq p\leq \infty$ % \leq q
$$
c \| \nabla ( \pi - \interp_T^{m-1} \pi)\|_{L^\infty(T)} \leq |T|^{-\frac 1 p} \| \nabla (\pi -\interp_T^{m-1}\pi) \|_{L^p(T)} \leq \| \nabla ( \pi - \interp_T^{m-1} \pi)\|_{L^\infty(T)}. %\leq  |T|^{-\frac 1 q} \| \nabla (\pi -\interp_T^{m-1}\pi) \|_{L^q(T)} 
$$
\end{lemma}
\begin{proof}
We denote by $T_0$ an arbitrary but fixed triangle of area $1$. For all $g \in L^\infty(T_0, \R^2)$ we obtain using Jensen's inequality
\be
\label{JensenLP}
\| g\|_{L^1(T_0)} \leq \| g \|_{L^p(T_0)} \leq \| g \|_{L^\infty(T_0)}.
\ee
Furthermore, since all norms are equivalent on the finite dimensional space $\P_{m-1}^2$, there exists a constant $c=c(m)>0$ such that
$
c \|\mu\|_{L^\infty(T_0)} \leq \| \mu \|_{L^1(T_0)}
$
 for all $\mu \in \P_{m-1}^2$. Therefore 
%\begin{gather}
\begin{IEEEeqnarray}{cCcCcCc}
\nonumber
c \|\mu\|_{L^\infty(T_0)} &\leq& \| \mu\|_{L^1(T_0)} &\leq& \| \mu \|_{L^p(T_0)} &\leq& \| \mu \|_{L^\infty(T_0)},\\
c \|\mu\|_{L^\infty(T)} &\leq& |T|^{-1}\| \mu\|_{L^1(T)} &\leq& |T|^{-\frac 1 p}\| \mu \|_{L^p(T)} &\leq& \| \mu \|_{L^\infty(T)},
\label{incNorms}
\end{IEEEeqnarray}
%$$
%\begin{array}{ccccccc} %{rrrrrrrrr}
%c \|\mu\|_{L^\infty(T_0)} &\leq& \| \mu\|_{L^1(T_0)} &\leq& \| \mu \|_{L^p(T_0)} &\leq& \| \mu \|_{L^\infty(T_0)},\\
%c \|\mu\|_{L^\infty(T)} &\leq& |T|^{-1}\| \mu\|_{L^1(T)} &\leq& |T|^{-\frac 1 p}\| \mu \|_{L^p(T)} &\leq& \| \mu \|_{L^\infty(T)},
%\end{array}
%$$
where we used in the second line a change of variables from $T_0$ to an arbitrary triangle $T$. We conclude the proof of this lemma by injecting $\mu = \nabla (\pi -\interp_T^{m-1}\pi)$ in \iref{incNorms}, where $\pi \in \H_m$ is arbitrary.
\end{proof}

For each $\mu \in \H_{m-1}^2$ we define 
\be
\label{defNormH2}
\|\mu \| := \sup_{z\neq 0} \frac{|\mu(z)|}{|z|^{m-1}},
\ee
and for each $A \in \M_2$ we define $\mu \circ A\in \H_{m-1}^2$ by $\mu \circ A(z) := \mu(A z)$. Note that $\|\mu\circ U \|= \|\mu\|$ for any $U \in \cO_2$, and that for $M \in S_2^+$
\be
\label{normMuM}
\|\mu \circ M^{-\frac 1 2}\| = \sup_{z\neq 0} \frac {|\mu(z)|}{|z|_M^{m-1}} \stext{ where } |z|_M := \sqrt{z^\trans M z},
\ee
hence for $M,M'\in S_2^+$
\be
\label{ineqMMpPi}
M \leq M' \ \Ra \ \| \pi\circ M^{-\frac 1 2 }\| \geq \| \pi\circ M'^{-\frac 1 2 }\|.
\ee
The next lemma shows that, among the triangles having a prescribed aspect ratio, encoded by a matrix $M$, the smallest approximation error is achieved by the \emph{acute triangles} (up to a fixed multiplicative constant) for which the measure of sliverness $S$ is minimal and equals one.

\begin{lemma}
\label{lemmaSM}
There exists a constant  $C=C(m)\geq 1$ such that the following holds. For any triangle $T$ and any $\pi \in \H_{m}$, 
one has denoting $M := \cH_T$
\be
\label{estimS}
C^{-1} \| (\nabla \pi) \circ M^{-\frac 1 2}\| \leq 
\|\nabla (\pi - \interp_T^{m-1} \pi)\|_{L^\infty(T)} \leq CS(T) \| (\nabla \pi) \circ M^{-\frac 1 2}\|.
\ee
\end{lemma}

\begin{proof}
We may assume that the triangle $T$ is centered at the origin of $\R^2$.
We obtain combining Lemma \ref{lemmaS} and the inclusion $T \subset \cE_T$, see \iref{defET},
$$
\|\nabla (\pi - \interp_T^{m-1} \pi)\|_{L^\infty(T)} \leq CS(T) \|\nabla \pi\|_{L^\infty(T)} \leq CS(T) \|\nabla \pi\|_{L^\infty(\cE_T)} = CS(T) \|  (\nabla \pi) \circ M^{-\frac 1 2}\|,
$$
which establishes the right part of \iref{estimS}.
On the other hand, if $C$ is sufficiently large, one has for all $\mu \in \H_{m-1}^2$
$$
C^{-1} \| \mu \| \leq \inf_{\nu \in \sP_{m-2}^2} \| \mu-\nu \|_{L^\infty(\TEq)},
$$
since both quantities are norms on $\H_{m-1}^2$. Since $z_T=0$, the change of variables $z\mapsto M^{-\frac 1 2} U^\trans z$ maps $\TEq$ onto $T$, where $U\in \cO_2$ is defined by \iref{eqMapTTeq}. Applying the above inequality to $\mu \circ (M^{-\frac 1 2}U^\trans)$, where $\mu\in \H_{m-1}^2$ is arbitrary, we thus obtain %obtain for any triangle $T$, denoting $M := \cH_T$
$$
C^{-1} \| \mu \circ M^{-\frac 1 2}\| = C^{-1} \| \mu \circ (M^{-\frac 1 2}U^\trans)\|  \leq \inf_{\nu \in \sP_{m-2}^2} \| \mu\circ (M^{-\frac 1 2}U^\trans)-\nu \|_{L^\infty(\TEq)}= \inf_{\nu \in \sP_{m-2}^2} \| \mu-\nu \|_{L^\infty(T)}.
$$
Choosing $\mu = \nabla \pi$ and $\nu = \nabla \interp_T^{m-1} \pi$ we obtain the left part of \iref{estimS}, which concludes the proof. % of this proposition.
\end{proof}

We introduce a variant $L_m$ of the shape function, following an idea originally proposed in \cite{C3}: for each $\pi\in \H_m$
\be
\label{defLm}
L_m(\pi) := \inf \{(\det M)^{\frac {m-1} 4} \ssep M\in S_2^+ \text{ s.t. } \|(\nabla \pi) \circ M\| \leq 1 \}.
\ee
\begin{lemma}
\label{lemmaLEquiv}
There exists a constant $C=C(m) \geq 1$ such that for all $\pi \in \H_m$
$$
C^{-1} L_m(\pi) \leq L_{m,p}(\pi) \leq CL_m(\pi).
$$
\end{lemma}
\begin{proof}
%This is immediately obtained by i
Injecting \iref{eqAreaHT}, Lemma \ref{lemmaEquiv} and Lemma \ref{lemmaSM} into \iref{optArea} we obtain
%\begin{eqnarray*}
%\inf \{ |T|^{-\frac {m-1} 2} \ssep T \text{ s.t.\ } |T|^{-\frac 1 p} \| \nabla (\pi-\interp_T^{m-1} \pi)\|_{L^p(T)} \leq 1\} &\leq& 
%\inf \{(\det M)^{\frac {m-1} 4} \ssep M\in S_2^+ \text{ s.t.\ } \|(\nabla \pi) \circ M\| \leq 1 \}\\
%&\leq&  \inf \{ |T|^{-\frac {m-1} 2} \ssep T \text{ s.t.\ } |T|^{-\frac 1 p} \| \nabla (\pi-\interp_T^{m-1} \pi)\|_{L^p(T)} \leq 1\},
%\end{eqnarray*}
\begin{IEEEeqnarray*}{CLlrl}
 &\inf \{ (|\TEq|/|T|)^{\frac {m-1} 2} &\ssep T \text{ s.t.\ } &|T|^{-\frac 1 p} &\| \nabla (\pi-\interp_T^{m-1} \pi)\|_{L^p(T)} \leq C\}\\
 \leq &
\inf \{(\det M)^{\frac {m-1} 4} &\IEEEeqnarraymulticol{2}{l}{\ssep M\in S_2^+ \text{ s.t.\ }} &\|(\nabla \pi) \circ M\| \leq 1 \}\\
\leq & \inf \{ (|\TEq|/|T|)^{\frac {m-1} 2} &\ssep T \text{ s.t.\ } & |T|^{-\frac 1 p}& \| \nabla (\pi-\interp_T^{m-1} \pi)\|_{L^p(T)} \leq c/C\},
\end{IEEEeqnarray*}
where we used the fact that for any $M\in S_2^+$ there exists an \emph{acute} triangle $T$ such that $\cH_T = M$.
Hence
$$
|\TEq|^{\frac {m-1} 2} L_{m,p}(\pi/C) \leq L_m(\pi) \leq |\TEq|^{\frac {m-1} 2} L_{m,p}(C \pi/c),
$$
which concludes the proof since the shape function satisfies the homogeneity property $L_{m,p}(\lambda \pi) = |\lambda |L_{m,p}(\pi)$ for any $\lambda\in \R$, as can be seen from the expression \iref{shapeFunctionL}.
\end{proof}

Following a suggestion of the referees, we may compare Lemma \ref{lemmaSM} to other anisotropic error estimates that have been proposed in the literature.
\begin{remark}
\label{remRefs}
The specificities of the estimate \iref{estimS} are the following:
\begin{enumerate}
\item It only applies to an \emph{exactly polynomial} function $f=\pi\in \H_m$, or $f=\pi+\mu$ where $\mu\in \P_{m-1}$ since the interpolation operator reproduces these elements: $\mu = \interp_T^{m-1} \mu$. This is a strong restriction, yet sufficient for our purposes, when compared for instance to \cite{Acosta} which applies to arbitrary functions in an adequate smoothness space.
\item It is \emph{sharp} (up to the multiplicative constant $C$) for any polynomial $\pi\in \H_m$ and any \emph{acute} triangle $T$. The sharpness is lost in the case of a strongly obtuse triangle since the measure of sliverness $S(T)$ only appears on the right of \iref{estimS}. A more elaborate estimate, which is restricted to the case $m=2$ of linear interpolation but applies to general functions, is proposed in \cite{ForPer01} in the attempt (confirmed by some examples) to provide a sharp estimate for both acute and obtuse triangles. A detailed study of the interpolation error of a polynomial function, with respect to the shape of the triangle $T$, can also be found in \cite{C4,C5} in the cases $m=2$ or $m=3$ respectively.
\item The orientation and the scales of the triangle $T$ are encoded in the matrix $M = \cH_T\in S_2^+$. This is one of the major strengths of this estimate, because it leads to the optimization problem \iref{defLm} posed on the set $S_2^+$ of symmetric positive definite matrices which can be addressed mathematically. This problem is studied in \S\ref{secShape} using algebraic techniques which yield an explicit equivalent of the shape functions $L_{m,p}$ for all $m \geq 2$, and an explicit ``near minimizer'' $\pi \in \H_m \mapsto \cM_m(\pi)\in S_2^+$, for $m\in\{2,3\}$, of the minimization problem \iref{defLm} defining $L_m(\pi)$. This problem is also studied in Chapter 6 of \cite{thesisJM} using analytical techniques, which yield well behaved, although implicit, ``near minimizers'' $\pi \mapsto \cM_m(\pi)$ of \iref{defLm} for all $m \geq 2$.
\end{enumerate}
\end{remark}

We now focus our attention on the global approximation of a function $f\in C^m(\overline \Omega)$ in the $W^{1,p}$ semi-norm. For that purpose we assume that a map $z\in \overline \Omega \mapsto \cM(z)\in S_2^+$ has been obtained which satisfies
\be
\label{defM}
\int_\Omega (\det \cM(z))^{\frac {m-1} 4\tau}dz \leq C_L^\tau \int_\Omega L_m(\pi_z)^\tau dz,  \stext{ and } \| \pi_z \circ \cM(z)^{-\frac 1 2}\| \leq 1 \text{ for all }  z\in \Omega,
\ee
where the polynomial $\pi_z\in \H_m$ is defined by \iref{taylorMuPi}, the exponent $\tau$ by \iref{defTau}, and $C_L$ is a constant not too large. In other words we assume that the matrix $\cM(z)$ is a minimizer, in an average sense and up to the sub-optimality constant $C_L$, of the optimization problem appearing in the definition \iref{defLm} of $L_m(\pi_z)$. Such a map can be obtained by setting $\cM(z) := \cM_m(\pi_z)+\delta \Id$, where $\cM_m$ is described in Point 3 of Remark \ref{remRefs}, and where $\delta>0$ is a positive constant introduced to avoid degeneracy problems.

\begin{figure}
\centering
\includegraphics[width=6cm,height=3.2cm]{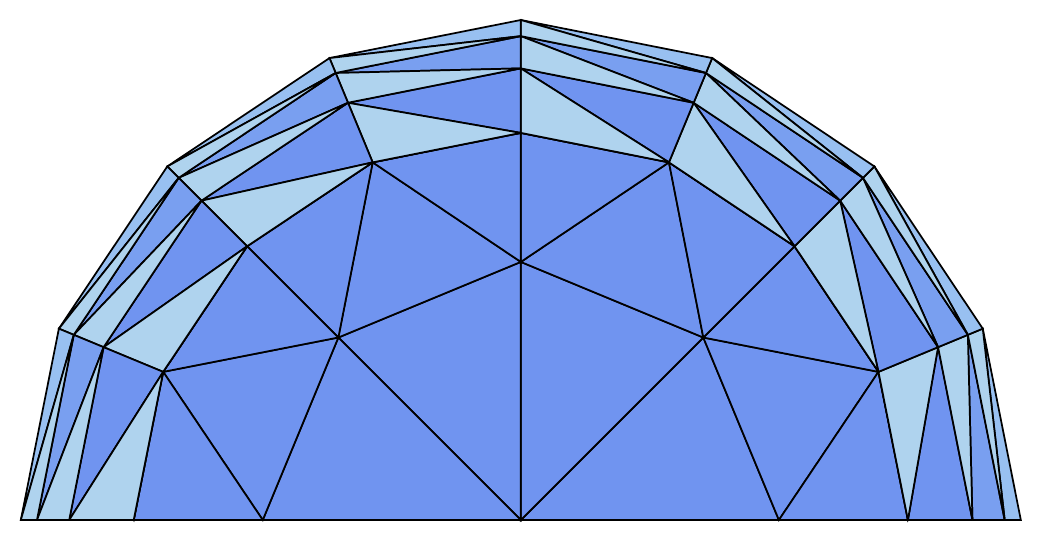}
\includegraphics[width=6cm,height=3.6cm]{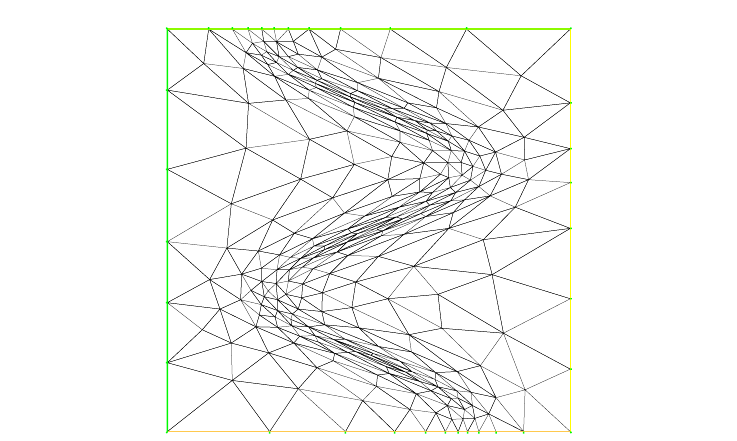}
\caption{Anisotropic mesh generation with a guaranteed upper bound on the maximal angle of the elements is generally tractable if one only requires some anisotropy \emph{close to the boundary} of the domain, and \emph{tangentially} to it (left). This is not any more the case if some anisotropy is required in the interior of the domain, at least with current software such as \cite{FreeFem} (right).}
\label{figMeshGen}
\end{figure}

We introduce a \emph{Riemannian metric} $H : \Omega \to S_2^+$
\be
\label{defH}
H(z) := h^{-2}(\det \cM(z))^\frac {-1}{(m-1)p+2} \cM(z),
\ee
where $h>0$ is a parameter. %, which eventually dictates the scale of the elements of the generated triangulation.
Mesh generation software such as \cite{FreeFem}, see Figure \ref{figMeshGen} (right), 
and \cite{Inria} in three dimensions (see \cite{Loseille} for a more extensive list), are designed to produce a mesh $\cT$ of $\Omega$ such that 
\be
\label{adaptTH}
C_0^{-2}H(z) \leq \cH_T \leq C_0^2 H(z)
\ee
for all $T\in \cT$, $z\in T$, where $C_0\geq 1$ is a constant not too large which reflects the quality of the adaptation of the mesh $\cT$ to the metric $H$. 
 In the expression \iref{defH} of the Riemannian metric $H$, the matrix $\cM(z)$ is used to prescribe optimal aspects ratios (requirement ii in the introduction) for the triangles $T\in \cT$, and the scalar factor to equidistribute the interpolation errors among the elements of $\cT$ (requirement i). Theoretical guarantees for such algorithms were established in \cite{Shew,Bois, thesisJM} when the metric sufficiently regular. %, and the property \iref{adaptTH} can be rigorously established under some conditions, see Chapter 5 of \cite{thesisJM}. 
Unfortunately these results \emph{do not} guarantee any property of the measure of sliverness $S(T)$ of the generated triangles $T\in \cT$ (requirement iii). This generally \emph{forbids} to achieve the optimal convergence estimate stated in Theorem \ref{mainTheorem}, even up to a fixed multiplicative constant. The adaptation of the method presented here to this (suboptimal) context is described \cite{sampTA}.\\

Limited results exist nevertheless on anisotropic mesh generation \emph{with some control} on the measure of sliverness.
Such a construction is (usually) possible if one only requires some anisotropy \emph{close to the boundary} of the domain, and \emph{tangentially} to it, see Figure \ref{figMeshGen} (left), which is the adequate behavior for the discretization of number of problems (such as the Poisson equation, or singularly perturbed reaction diffusion problems) as discussed in \cite{ApKu}.
In the general case where some anisotropy is required in the interior of the domain, we refer to Theorem 6.1.2 in Chapter 6 of \cite{thesisJM}, on ``quasi-acute triangulations'', which implies the following.
Assume that $\Omega$ is the periodic domain $(\R/\Z)^2$; the key is the absence of a boundary. Assume that  the metric $H$ satisfies the strong regularity properties required for this result, which are expressed under the form of Lipschitz regularity with respect to some distances on $S_2^+$ and $\R^2$, and are satisfied if the parameter $h>0$ in \iref{defH} is sufficiently small according to Lemma 6.5.12 in \cite{thesisJM}. Then there exists a triangulation $\cT$ which satisfies \iref{adaptTH}, for an absolute constant $C_0$ independent of $H$, and in addition the following property: there exists a $C_0$-refinement $\cT'$ of the triangulation $\cT$ (in other words each element of $\cT$ contains at most $C_0$ elements of $\cT'$, and each element of $\cT'$ is contained in an element of $\cT$) such that $S(T') \leq C_0$ for all $T'\in \cT'$ (requirement iii in the introduction). 

Let $(\cT, \cT')$ be a pair of triangulations satisfying the above conditions.
If $T(z)\in \cT$ denotes the element containing a point $z\in \Omega$, then we obtain using \iref{eqAreaHT} and \iref{adaptTH} 
\begin{IEEEeqnarray}{cCcCcCcCc}
\#(\cT') &\leq& C_0\#(\cT) &=& C_0 \int_\Omega \frac {dz}{|T(z)|} &\leq& \frac{C_0^3}{|\TEq|} \int_\Omega \sqrt {\det H} &=& \frac{C_0^3}{|\TEq|} h^{-2} \int_\Omega (\det \cM)^{\frac {m-1} 4\tau}.
\label{cardTp}
%\#(\cT') \leq C_0\#(\cT) = C_0 \int_\Omega \frac {dz}{|T(z)|}\leq \frac{C_0^3}{|\TEq|} \int_\Omega \sqrt {\det H} = \frac{C_0^3}{|\TEq|} h^{-2} \int_\Omega (\det \cM)^{\frac {m-1} 4\tau}.
\end{IEEEeqnarray}
Let $z\in \Omega$, $T'\in \cT'$ and $T\in \cT$ be such that $z\in T'\subset T$. Using \iref{ineqTH} and \iref{adaptTH} we obtain $\cH_T' \geq 2^{-2} \cH_T \geq (2 C_0)^{-2}H(z)$, hence 
%For each  $T' \in \cT'$ and each $z\in T'$ one obtains $\cH_T' \geq (2 C_0)^{-2}H(z)$, using \iref{ineqTH} and \iref{adaptTH}, hence
%Furthermore one has for each $T \in \cT'$ and each $z\in T$ one has, observing that $H(z) \leq (2 C_0)^2\cH_T'$
\begin{eqnarray*}
%\|\nabla (f-\interp_T^{m-1} f)\|_{L^p(T)} &\lsim& 
\| \nabla (\pi_z - \interp_{T'}^{m-1} \pi_z)\|_{L^p(T')}
 &\leq&  |T'|^\frac 1 p \| \nabla (\pi_z - \interp_{T'}^{m-1} \pi_z)\|_{L^\infty(T')}\\
 &\leq& C |\TEq|^\frac 1 p (\det \cH_{T'})^{-\frac 1 {2p}} S(T') \| (\nabla \pi_z) \circ \cH_{T'}^{-\frac 1 2}\|\\
 &\leq & C_1 (\det H(z))^{-\frac 1 {2p}} \| (\nabla \pi_z) \circ H(z)^{-\frac 1 2}\|\\
 &=& C_1 h^{\frac 2 \tau} \|(\nabla \pi_z) \circ \cM(z)^{-\frac 1 2}\|\\
 &\leq& C_1 h^{\frac 2 \tau},
% &\leq & C |\TEq|^\frac 1 p (2C_0)^{\frac 1 \tau} S(T')  (\det H(z))^{-\frac 1 {2p}} \| (\nabla \pi_z) \circ H(z)^{-\frac 1 2}\|\\
% &=& C |\TEq|^\frac 1 p (2C_0)^{\frac 1 \tau} S(T') h^{\frac 2 \tau} \|(\nabla \pi_z) \circ \cM(z)^{-\frac 1 2}\|\\
% &\leq& C |\TEq|^\frac 1 p 2^\frac 1 \tau C_0^{\frac 1 \tau+1} h^{\frac 2 \tau},
\end{eqnarray*}
where $C_1 := C |\TEq|^\frac 1 p 2^\frac 1 \tau C_0^{\frac 1 \tau+1}$. We used successively Jensen's inequality in the first line, \iref{eqAreaHT} and Lemma \ref{lemmaSM} in the second line, $S(T') \leq C_0$ the definition \iref{defTau} of $\tau$ and \iref{ineqMMpPi} in the third line, \iref{defH} in the fourth line and \iref{defM} in the last line.
In order to control the function $f$, instead of the polynomial $\pi_z$ which corresponds to its Taylor expansion, we may use Point i of Lemma \ref{lemmaDiff}, proven in \S\ref{secProofs} below, which immediately implies in this context that 
$$
|\|\nabla (f-\interp_{T'}^{m-1}f)\|_{L^p(T')} - \| \nabla (\pi_z - \interp_{T'}^{m-1} \pi_z)\|_{L^p(T')}| \leq \ve(h) h^\frac 2 \tau,
$$
where $\ve(h) \to 0$ as $h \to 0$, and the function $\ve : \R_+^*\to \R_+$ only depends on $f$, $m$ and $\cM$.
If $h$ is sufficiently small (requirement iv in the introduction), then $\ve(h) \leq 1$ and therefore $\|\nabla (f-\interp_{T'}^{m-1}f)\|_{L^p(T')}\leq (C_1+1) h^\frac 2 \tau$ for each $T'\in \cT'$. It follows that 
\begin{eqnarray*}
\#(\cT')^{\frac {m-1} 2} \|\nabla (f-\interp_{\cT'}^{m-1} f)\|_{L^p(\Omega)} &\leq& \#(\cT')^{\frac 1 \tau} \max_{T'\in \cT'} \|\nabla (f-\interp_{T'}^{m-1} f)\|_{L^p(T')}\\
&\leq& (C_1+1) (C_0^3/|\TEq|)^\frac 1 \tau \left(\int_\Omega (\det \cM)^{\frac {m-1} 4 \tau} dz\right)^\frac 1 \tau,\\
& \leq & C_2 \|L_m(d^m f/m!)\|_{L^\tau(\Omega)},
\end{eqnarray*}
where $C_2 = (C_1+1) (C_0^3/|\TEq|)^\frac 1 \tau C_L$.
We used the definition \iref{defTau} of the exponent $\tau$ in the first line, \iref{cardTp} in the second line and \iref{defM} in the last.
As announced the triangulation $\cT'$ satisfies the optimal error estimate of Theorem \ref{mainTheorem}, up to the multiplicative constant $C_2$. 

%Yet another anisotropic mesh generation method is the \emph{hierarchical} refinement procedure presented in \cite{CDHM}, and proved to be optimal in \cite{CMi} in the special case of the piecewise linear approximation of a convex function in the $L^p$ norm. 
%%%This method produces 
%%The sequences of anisotropic triangulations $(\cT_N)_{N \geq 0}$ produced with this method are hierarchical in the sense that $\cT_{N+1}$ is a refinement $\cT_N$ for each $N \geq 0$, a property which is desirable in a number of applications and allows as a by-product to define an anisotropic wavelet basis. 
%%Unfortunately, again, the elements of these meshes may have a largest angle arbitrarily close to $\cPi$.
%
%\begin{remark}
%The mesh $\cT$ of the domain $\Omega$ can be designed in a \emph{structured} way in some situations. This is typically the case is the domain $\Omega$ is simple and if anisotropy is only required close to the boundary $\partial \Omega$, and tangentially to it. Ensuring that the angles are bounded away from $\cPi$ is generally possible in this situation, as illustrated on Figure \ref{figMeshGen} (left). Such constructions, which also apply in three dimensions, are particularly well tailored to Laplace and Maxwell equations on singular domains, see \cite{??}. %, for which exhibit this type of anisotropic behavior and for which 
%\end{remark}

\section{Study of the shape function}
\label{secShape}
This section is devoted to the close study of the shape function $L_{m,p}$, using algebraic techniques. %: we obtain explicit \emph{near} minimizers  We obtain explicit Our starting point is the variant %, or more the variant introduced in \iref{defLm}, defined for all $\pi \in \H_m$ by 
Our approach is based on the variant $L_m$ introduced in \iref{defLm}, defined for all $\pi \in \H_m$ by 
%the variant $L_m$ introduced in \iref{defLm} and defined for all $\pi \in \H_m$ by
\be
\label{defLm2}
L_m(\pi) := \inf\{(\det M)^{\frac {m-1} 4} \ssep M \in S_2^+ \text{ s.t. } \|(\nabla \pi) \circ M^{-\frac 1 2} \| \leq 1\},
\ee
and which is uniformly equivalent to $L_{m,p}$ according to Lemma \ref{lemmaLEquiv}.

\subsection{Explicit minimizers}
%We describe  a close to optimal triangle for the interpolation of each polynomial $\pi\in \H_m$, when $m \in \{2,3\}$: among all triangles of area $1$, this triangle minimizes the approximation error $\|\nabla (\pi- \interp_T^{m-1})\|_{L^p(T)}$ up to a fixed multiplicative constant.
%
%
%If $M\in S_2^+$ is a minimizer of the optimization problem appearing in \iref{defLm} up to the multiplicative constant $C_0$, then any acute triangle $T$ satisfying $\cH_T = M$ is a minimizer of the optimization problem appearing in  \iref{optArea} up to the multiplicative constant $C_0C_1$, where $C_1 = C_1(m)$ is independent of $\pi$.

We describe the solution to the optimization problem appearing in \iref{defLm2} when $m\in \{2,3\}$. The case $m=2$, which corresponds to piecewise linear finite elements is already known, and discussed in detail in \cite{Shew,AzSimp91} and \cite{HuSu03} for functions of more than two variables. 
In contrast the results in the case $m=3$ are entirely new, although this case, which corresponds to piecewise quadratic elements, has already been discussed in \cite{C5}.
%The case $m=3$, which corresponds to piecewise quadratic elements, has already been discussed in \cite{C5}, but the results presented below are entirely new. 

In order to present our results, we introduce some notation.
For any homogeneous quadratic polynomial $\pi\in \H_2$, $\pi = a x^2 + 2 b x y + c y^2$, we define the symmetric matrix
\be
\label{defMatPi}
 [\pi] = 
\left(\begin{array}{cc}
a & b\\
b & c
\end{array}\right).
\ee
For all $\pi \in \H_2$,  $\pi = a x^2 + 2 b x y + c y^2$, we define
\be
\label{defMP1}
\cM_2(\pi) := 4 [\pi]^2 = 
4 \left(\begin{array}{cc}
a & b\\
b & c
\end{array}\right)^2
= 
4 \left(\begin{array}{cc}
a^2+b^2 & ab+bc\\
ab+bc & b^2+c^2
\end{array}\right).
\ee
For all $\pi\in \H_3$, $\pi = a x^3+ 3 b x^2 y+ 3 c x y^2 + d y^3$, we define
\be
\label{defMP2}
\cM_3(\pi) := \sqrt{[\partial_x \pi]^2+[\partial_y \pi]^2} = 
3 \sqrt{
\left(\begin{array}{cc}
a & b\\
b & c
\end{array}\right)^2
+
\left(\begin{array}{cc}
b & c\\
c & d
\end{array}\right)^2
}
= 3
\sqrt{
\left(\begin{array}{cc}
a^2+2b^2+ c^2 & ab+2 bc + cd\\
ab+2 bc + cd & b^2+2c^2+d^2
\end{array}\right).
}
\ee
We say that a polynomial $\pi \in \H_m$ is \emph{univariate} if there exists $\lambda, \alpha, \beta\in \R$ such that $\pi = \lambda (\alpha x+ \beta y)^m$.
\begin{prop}
\label{propcM}
\begin{enumerate}[i.]
\item If $\pi \in \H_2$ is not univariate, then the matrix $\cM_2(\pi)$ is the unique minimizer of the optimization problem appearing in \iref{defLm}.
\item The map $\pi\in \H_3 \to \cM_3(\pi)$ is a \emph{near-minimizer} of the problem \iref{defLm} in the following sense. If $\pi\in \H_3$ is not univariate, then $\cM_3(\pi)$ is non-degenerate and $\| (\nabla \pi) \circ \cM_3(\pi)^{-\frac 1 2}\| \leq \sqrt 2$ (hence $\| (\nabla \pi) \circ (\sqrt 2\cM_3(\pi))^{-\frac 1 2}\| \leq 1$). Furthermore there exists a constant $C$, independent of $\pi$, such that
\be
\label{minM3}
\sqrt{\det \cM_3(\pi)} \leq C L_3(\pi). %\inf\{ \det M \sep M\in S_2^+ \text{ and } \| (\nabla \pi) \circ M^{-\frac 1 2}\| \leq 1\}.
\ee
\end{enumerate}
\end{prop}

\begin{proof}
According to \iref{normMuM}, $\| (\nabla \pi) \circ M^{-\frac 1 2}\| \leq K$ is equivalent to $|\nabla \pi(z)| \leq K|z|_M^{m-1}$ for all $z\in \R^2$,  where $\pi \in \H_m$, $M\in S_2^+$ and $K>0$ are arbitrary.

We first consider an homogeneous polynomial $\pi \in \H_2$, which is not univariate.
For all $z\in \R^2$ one has $\nabla \pi(z) = 2 [\pi] z$, and therefore $|\nabla \pi(z)| ^2 = z^\trans \cM_2(\pi) z$. On the other hand consider $M\in S_2^+$ such that $|\nabla \pi(z)|^2 = z^\trans \cM_2(\pi) z\leq z^\trans M z$. We thus have  $\cM_2(\pi) \leq M$ in the sense of symmetric matrices, which implies that $\det \cM_2(\pi) \leq \det M$ with equality if and only if $\cM_2(\pi) = M$, since $\cM_2(\pi)$ is positive definite. This concludes the proof of the first point.

We now consider $\pi \in \H_3$, which is again not univariate. 
In the sense of symmetric matrices, we have
$$
\cM_3(\pi) = \sqrt{[\partial_x \pi]^2+[\partial_y \pi]^2} \geq \sqrt{[\partial_x \pi]^2} = |[\partial_x \pi]|,
$$
where we used the fact that the square root $\sqrt\cdot : S_2^\oplus\to S_2^\oplus$ is increasing.
It follows that 
$$
|\nabla \pi(z)| ^2 =  |\partial_x\pi(z)|^2 + |\partial_y\pi(z)|^2 \leq 2 (z^\trans \cM_3(\pi) z)^2,
$$
hence $\cM_3(\pi)$ satisfies the constraint $\| (\nabla \pi) \circ \cM_3(\pi)^{-\frac 1 2}\| \leq \sqrt 2$. 
Note that
\be
\label{detMPi}
\det \cM_3(\pi) = 9 \sqrt{(a^2+2b^2+c^2)(b^2+2c^2+d^2) - (ab+2bc+cd)^2}.
\ee
We postpone the proof of \iref{minM3} to \S \ref{secPolEq}, right after \iref{L3eq}, as we develop a general method for obtaining simple equivalents of the functions $L_m$. 
\end{proof}

%\begin{remark}
%Let $\pi \in \H_m$, be non univariate, let $\cM\in S_2^+$ be such that 
%Produce $T$ of area $1$ from $\cM_m$. (*details*)
%%If $M\in S_2^+$ is a minimizer of the optimization problem appearing in \iref{defLm} up to the multiplicative constant $C_0$, then according to Lemma \ref{lemmaSM} any acute triangle $T$ satisfying $\cH_T = M$ is a minimizer of the optimization problem appearing in  \iref{optArea} up to the multiplicative constant $C_1$, where $C_1 = C_1(m)$ is independent of $\pi$. It then follows from \iref{transInv} that the triangle $T'$ homothetic to $T$ and of area $1$ satis
%\end{remark}

\begin{remark}
It was proposed in \cite{AzSimp91,HuSu03} to generate an anisotropic mesh $\cT$ of a domain $\Omega$ via the transport of a \emph{uniform} mesh $\cT'$ of an auxiliary domain $\Omega'$ by a diffeomorphism $F : \Omega \to \Omega'$. Without entering the details of this approach, we may describe one of its successes. 
%The following mesh generation method, the transport of a uniform mesh of a domain $\Omega'$ onto the domain of interest $\Omega$ via a diffeomorphism, has been presented (slightly differently) and generalized in \cite{AzSimp91,HuSu03}. %, presented in a slightly different manner.
Assume that one wishes to approximate a strongly convex function $f\in C^2(\overline \Omega)$, using linear finite elements. Define $F(z) := 2\nabla f(z)$ and $\Omega' := F(\Omega) \subset \R^2$. Consider a \emph{uniform} mesh $\cT'$ of $\Omega'$, of mesh size $h>0$, and denote by $\cT$ the collection of triangles 
%Mesh displacement is another anisotropic mesh generation method, see [??],which proceeds as follows: consider an auxiliary domain $\Omega'$, a \emph{uniform} triangulation $\cT'$ of $\Omega'$, and a one to one and $C^1$ map $F : \Omega \to \Omega'$. Denote by $\cT$ the collection of triangles 
obtained as follows: for each triangle $T'\in \cT'$, of vertices $v_1,v_2,v_3$, the set $\cT$ contains the triangle $T$ of vertices $F^{-1}(v_1), F^{-1}(v_2), F^{-1}(v_3)$. If the parameter $h>0$ is sufficiently small, then $\cT$ is a triangulation since $F: \Omega \to \Omega'$ is a diffeomorphism. Furthermore we obtain using \iref{eqInvHT} that $\cT$ is adapted in the sense of \iref{adaptTH} to the metric $H : \Omega \to S_2^+$ defined by
%If the function $F$ is well behaved then $\cT$ is a triangulation of $\Omega$, which is equivalent in the sense of \iref{adaptTH} to the metric 
\be
\label{eqHF}
H(z) := h^{-2} dF(z)^\trans dF(z) = h^{-2} 4 (d^2 f(z))^2 = h^{-2}\cM_2(\pi_z),
\ee
hence the \emph{optimal metric}, see \iref{defH}, for the approximation $f$ in the $W^{1,\infty}$ semi-norm! (The control of the measure of sliverness of the elements of $\cT$ remains, however, an open problem with this method.)%Unfortunately the measure of sliverness of the elements of $\cT$ is hard to control. %Furthermore this method does not allow to produce a triangulation adapted in the sense of \iref{adaptTH} to an \emph{arbitrary} Riemannian metric: the Riemann-Christoffel symbol, a differential invariant of the metric, needs to vanish identically.
%where $h>0$ is a parameter tied to the density of the mesh $\cT'$ of $\Omega'$. This observation results from the transformation rule \iref{eqInvHT} for the matrix $\cH_T$, when a linear change of coordinates is applied to the triangle $T$. This method not universal : for a given metric $H$ on $\Omega$, there exists a domain $\Omega'$ and a map $F:\Omega\to \Omega'$ satisfying \iref{eqHF} only if the Riemann-Christoffel symbol of $H$ vanishes uniformly on $\Omega$. Yet it is successful, as observed in \cite{??}, in the case of the piecewise linear approximation in the $W^{1,\infty}$ semi-norm of a $C^2$ strongly convex function $f$. For which the optimal metric $H$ is $h^{-2}\cM_2(\pi_z) = h^{-2} (d^2f(z))^2$ has vanishing Riemann-Christoffel symbols, and the corresponding map is $F = \nabla f$. The author has not heard so far of guaranteed upper bound on the maximal angle of the elements generated with this method.
\end{remark}

Let us finally mention that, although they are derived from the
coefficients of $\pi$, the maps $\pi \mapsto \cM_m(\pi)$ 
for $m\in\{2,3\}$ are invariant under rotation,
and therefore not tied to the chosen system of
coordinate $(x,y)$, as expressed by the following result.

\begin{prop}
For any $m\in\{2,3\}$, any $\pi\in \H_m$ and any unitary matrix $U\in \cO_2$, one has
$$
\cM_m(\pi\circ U) = U^\trans \cM_m(\pi) U.
$$
\end{prop}

\begin{proof}
Let $\pi\in \H_2$ and let $U\in \cO_2$. Then, as announced,
$$
\cM_2(\pi\circ U) = 4[\pi\circ U]^2 = 4 (U^\trans [\pi]U)^2 = U^\trans \cM_2(\pi) U.
$$
Let $\pi\in \H_3$ and denote by $(u_{ij})_{1 \leq i,j\leq 2}$ the entries of the unitary matrix $U$. Then 
%let 
%$
%U = 
%\left(
%\begin{array}{cc}
%u_{11} & u_{12}\\
%u_{21} & u_{22}
%\end{array}
%\right)
%$
%be unitary,
%then 
$$
[\partial_x(\pi\circ U)] = u_{11} U^\trans [\partial_x \pi] U + u_{12} U^\trans [\partial_y \pi] U \quad \text{ and } \quad [\partial_y(\pi\circ U)] = u_{21} U^\trans [\partial_x \pi] U + u_{22} U^\trans [\partial_y \pi] U
$$
Hence, since $U$ is unitary,
\begin{align*}
&[\partial_x(\pi\circ U)]^2 +  [\partial_y(\pi\circ U)]^2\\
& = (u_{11}^2+u_{21}^2) U^\trans [\partial_x \pi]^2 U + (u_{11} u_{12}+u_{21} u_{22}) U^\trans ([\partial_x \pi] [\partial_y \pi]+ [\partial_y \pi] [\partial_x \pi]) U 
+ (u_{12}^2+u_{22}^2) U^\trans [\partial_y \pi]^2 U\\
& = U^\trans [\partial_x \pi]^2 U+U^\trans [\partial_y \pi]^2 U
\end{align*}
%$$
%[\partial_x(\pi\circ U)]^2 +  [\partial_y(\pi\circ U)]^2 = (u_{11}^2+u_{21}^2) U^\trans [\partial_x \pi]^2 U + (u_{11} u_{12}+u_{21} u_{22}) U^\trans ([\partial_x \pi] [\partial_y \pi]+ [\partial_y \pi] [\partial_x \pi]) U 
%+ (u_{12}^2+u_{22}^2) U^\trans [\partial_y \pi]^2 U
%$$
%which equals $U^\trans [\partial_x \pi]^2 U+U^\trans [\partial_y \pi]^2 U$  
Therefore %Eventually
$$
\cM_3(\pi\circ U) = \sqrt{U^\trans [\partial_x \pi]^2 U+U^\trans [\partial_y \pi]^2 U} = U^\trans \sqrt{[\partial_x \pi]^2+[\partial_y \pi]^2} \,U = U^\trans \cM_3(\pi) U
$$
which concludes the proof.
\end{proof}

\subsection{Polynomial equivalents}
\label{secPolEq}
%The optimal error estimates established in Theorem \ref{mainTheorem}
%involve the %quantity $L_{m,p}(\frac {d^m f}{m!})$. The
%shape function 
%$\pi \mapsto L_{m,p}(\pi)$ which is obtained by solving an optimization problem,
%and does not have an explicit analytic expression 
%in terms of the coefficients of $\pi\in \H_m$.
We introduce equivalents of the shape function $\pi\mapsto L_m(\pi)$ on $\H_m$, %we introduce quantities which
%are equivalent to $L_m(\pi)$% defined by \iref{defLm}, and therefore to $L_{m,p}(\pi)$ for all $1\leq p\leq \infty$ according to Lemma \ref{lemmaLEquiv}, 
 %and  
 which can be written in analytic form in terms of the coefficients of $\pi\in \H_m$.
As a starter we infer from Point i of Proposition \ref{propcM} that for any $\pi\in \H_2$ %, that for each $\pi\in \H_2$
$$
L_2(\pi) = (\det \cM_2(\pi))^\frac 1 4 = 2\sqrt {|\det [\pi]|}.
$$
%where the matrices $[\pi]\in S_2$ and $\cM_2(\pi)\in S_2^\oplus$ are defined respectively by \iref{defMatPi} and \iref{defMP1}.

For each $r\geq 2$ we denote by $\H_r$ the space of homogeneous bivariate polynomials of degree $r$, as in \iref{defPH} for $r=m$, equipped with the norm 
$$
\|\mu\| := \sup_{z\neq 0} \frac {|\mu(z)|}{|z|^r}.
$$
For $\mu\in \H_r$ and $A \in \M_2$, we define the homogeneous polynomial $\mu\circ A\in \H_r$ by $\mu\circ A(z) := \mu(Az)$, $z\in \R^2$. Observe that for $M\in S_2^+$ one has
$$
\|\mu \circ M^{-\frac 1 2}\| = \sup_{z\neq 0} \frac {|\mu(z)|}{|z|_M^r}.
$$
We now introduce a variant $K_r := \H_r \to \R_+$ of the shape function $L_m$, which was first defined in \cite{C3}, and later studied in \cite{Mi}, in the context of optimal mesh adaptation for finite element approximation in the $L^p$ norm. %, when the error is measured in the $L^p$ norm instead of the $W^{1,p}$ semi-norm considered here.  
(More precisely, due to different conventions, the function $K_r$ is tied to the function $K_r^\cE$ defined in \cite{Mi} by the equality $K_r = \cPi^{-\frac r 2} K_r^\cE$.) For each $\mu\in \H_r$ we define
$$
K_r(\mu) := \inf\{ (\det M)^\frac r 4\ssep M \in S_2^+ \text{ s.t. }\|\mu \circ M^{-\frac 1 2}\| \leq 1\}.
$$
Observe that $|\nabla \pi|^2 := (\partial_x \pi)^2+ (\partial_y \pi)^2 \in \H_{2m-2}$ for each $\pi \in \H_m$, and that clearly 
\be
\label{eqLK}
L_m(\pi) = \sqrt{K_{2m-2}(|\nabla \pi|^2)}.
\ee

Given a pair of non negative functions $Q$ and $R$ on $\H_m$, we write $Q\sim R$ if and only if there exists a constant $C\geq 1$ such that $C^{-1} Q\leq R \leq C Q$ holds uniformly on $\H_m$. 
We sometimes slightly abuse notations and write $Q(\pi)\sim R(\pi)$.
We say that a function $Q$ is a polynomial on $\H_m$ if there exists a polynomial $P$ of $m+1$ real variables such that for all $a_0, \cdots,  a_m\in \R$,
$$
Q\left(\sum_{i=0}^m a_i x^i y^{m-i}\right) = P(a_0,\cdots, a_m).
$$
We define $\deg Q := \deg P$, and we say that $Q$ is homogeneous
if $P$ is homogeneous. 

The following equivalences were established in \cite{Mi} : for $\pi \in \H_2$
\be
\label{K2}
K_2(\pi) \sim \sqrt{|\det [\pi]|},
\ee
and for $\pi \in \H_3$
\be
\label{K3}
K_3(\pi)\sim \sqrt[4]{|\disc(\pi)|},
\ee
where $\disc(\pi)$ denotes the discriminant of the cubic polynomial $\pi$, which is defined by 
$$
\disc(a x^3 + 3b x^2  y+ 3c x y^2 + dy^3) =  4 (ab-c^2)(bc-d^2) - (ad-bc)^2. %b^2c^2 - 4ac^3 - 4b^3d + 18abcd - 27a^2 d^2.
$$
More generally for each $r\geq 2$ an explicit homogenous polynomial $Q$ on $\H_r$ is introduced in \cite{Mi}, which satisfies 
$$
K_r \sim \sqrt[d]{|Q|}, \text{ with } d := \deg Q.  
$$
Combining this result with \iref{eqLK} we obtain an explicit equivalent of the functions $L_m$.
\begin{prop}
\label{thPolEq}
Let $m\geq 2$ and let $Q$ be an homogeneous polynomial on $\H_{2m-2}$ such that 
$K_{2m-2}\sim \sqrt[d]{|Q|}$, where $d = \deg Q$. 
Let $Q_*$ be the polynomial on $\H_m$ defined by %defined for all $\pi\in\H_m$ by 
$$
Q_*(\pi) := Q(|\nabla \pi|^2),
$$
then $L_m \sim \sqrt[2d]{Q_*}$ on $\H_m$.
\end{prop}

This construction uses an equivalent $\sqrt[d]{|Q|}$ of $K_{2m-2}$ to produce an equivalent of $L_m$.
Unfortunately, as $m$ increases, the practical construction of $Q$ becomes more involved and the degree
$d$ quickly raises. In the following theorem, we build 
an equivalent to $L_{m}$ from an equivalent of $K_{m-1}$ instead of $K_{2m-2}$, which is therefore
simpler.

\begin{theorem}
\label{thPolEq2}
Let $m\geq 3$ and let $Q$ be an homogeneous polynomial on $\H_{m-1}$ such that $K_{m-1} \sim \sqrt[d]{|Q|}$, 
where $d=\deg Q$. Let $(Q_k)_{0\leq k\leq r}$ be the 
homogeneous polynomials of degree $d$ on $\H_{m-1}\times \H_{m-1}$ such that 
for all $u,v\in \R$ and all $\pi_1,\pi_2\in \H_m$ one has
\be
\label{defQk}
Q(u \pi_1+v \pi_2) = \sum_{0\leq k\leq d} \binom d k u^k v^{d-k} Q_k(\pi_1,\pi_2),
\ee
where $\binom d k := \frac {d!}{k!(d-k)!}$.
Let $Q_*$ be the polynomial defined for all $\pi\in\H_m$ by 
$$
Q_*(\pi) := \sum_{0\leq k\leq d} \binom d k Q_k\left(\partial_x \pi,\partial_y \pi\right)^2,
$$ 
then $L_{m} \sim \sqrt[2d]{Q_*}$ on $\H_{m}$.
\end{theorem}

\begin{proof}
See Appendix \ref{appenLQ}.
\end{proof}

Using this construction and \iref{K2} we obtain an equivalent of $L_3$ as follows.
Let $\pi_1 = a x^2+ 2b x y + c y^2$ and $\pi_2 = a' x^2+ 2b' x y + c' y^2$ be two elements
of $\H_2$. We obtain
\begin{align*}
\det( [u \pi_1 + v \pi_2]) &= (ua+va') ( uc+vc') - (ub+vb')^2\\
& = u^2 (ac-b^2) + u v (ac'+a'c - 2 b b') + v^2 (a'c' - b'^2).
\end{align*}
Applying the construction of Theorem \ref{thPolEq2} to $\pi = a x^3 + 3 b x^2  y+ 3 c x y^2 + dy^3\in \H_3$ we obtain
\be
\label{L3eq}
L_3(\pi) \sim 3 \sqrt[4]{(ac-b^2)^2+(ad-bc)^2/2+(bd-c^2)^2}.
\ee
Remarking that 
$$
2 [(ac-b^2)^2+(ad-bc)^2/2+(bd-c^2)^2 ]= (a^2+2b^2+c^2)(b^2+2c^2+d^2) - (ab+2bc+cd)^2,
$$
 and using equation \iref{detMPi} we obtain that $L_3(\pi) \sim \sqrt{\det \cM_3(\pi)}$. 
This point concludes the proof of Proposition \ref{propcM} and thus establishes that the map $\cM_3$ defined in \iref{defMP2} can be used for optimal mesh adaptation for quadratic finite elements.
\newline
\newline
\noindent
Using \iref{K3} we similarly 
obtain an equivalent of $L_4$. Denoting $\pi = a x^4 + 4 bx^3 y+ 6c x^2 y^2+ 4 d x y^3 + e y^4$: 
\begin{align*}
L_4(\pi)^8  \sim {} & (3 b^2 c^2 - 4 a c^3 - 4 b^3 d + 6 a b c d - a^2 d^2)^2 \\
&+ (2 b c^3 - 6 a c^2 d + 4 a b d^2 - 4 b^3 e + 6 a b c e - 2 a^2 d e)^2/4 \\
&+ (3 c^4 - 6 b c^2 d + 8 b^2 d^2 - 6 a c d^2 - 6 b^2 c e + 6 a c^2 e + 2 a b d e - a^2 e^2)^2/6\\
&+ (2 c^3 d - 4 a d^3 - 6 b c^2 e + 4 b^2 d e + 6 a c d e - 2 a b e^2)^2/4 \\
&+ (3 c^2 d^2 - 4 b d^3 - 4 c^3 e + 6 b c d e - b^2 e^2)^2.
\end{align*}

The following proposition identifies the polynomials $\pi\in \H_m$ for which $L_m(\pi) = 0$, and therefore the 
values of $d^{m}f$ for which anisotropic mesh adaptation may lead to \emph{super-convergence}.

\begin{prop}
Let $m\geq 2$ and let $t_m := \left\lfloor \frac{m+3} 2 \right\rfloor$. Then for all $\pi \in \H_m$,
\be
\label{vanishL}
L_m(\pi) = 0 \text{ if and only if } \pi = (\alpha x+ \beta y)^{t_m} \tilde \pi \text{ for some } \alpha, \beta\in \R \text{ and } \tilde \pi \in \H_{m-t_m}.
\ee
\end{prop}
\proof
It was established in \cite{Mi} that $K_{2m-2}(\pi_*)=0$ if and only if $\pi_*\in \H_{2m-2}$ has a linear factor of multiplicity $m$. We therefore obtain, using \iref{eqLK}, that  $L_m(\pi) = 0$ if and only
$|\nabla \pi|^2$ is a multiple of $l^m$, where $l$ is of the form $l=\alpha x +\beta y$.

Let us first assume that $|\nabla \pi|^2 = (\partial_x \pi)^2+ (\partial_y \pi)^2$ has such a form. Clearly  $(\partial_x \pi)^2$ and $(\partial_y \pi)^2$ are both multiples of $l^m$. Therefore $\partial_x \pi$ and $\partial_y \pi$ are multiples of $l^s$ where $s$ is an integer such that $2s\geq m$, hence $s\geq t_m-1$. We therefore have 
$$
\partial_x \pi = l^s \pi_1 \ \text{ and } \ \partial_y \pi = l^s \pi_2 \text{ where } \pi_1,\pi_2 \in \H_{m-1-s}.
$$
Recalling that $l = \alpha x+\beta y$ we obtain
$$
0 = \partial_{yx}^2 \pi -\partial_{xy}^2 \pi = l^s(\partial_y \pi_1 - \partial_x\pi_2) + sl^{s-1} (\beta\pi_1-\alpha\pi_2),
$$
hence $\beta\pi_1-\alpha\pi_2$ is a multiple of $l$. Since $\pi$ is homogenous of degree $m$ it obeys
the Euler identity $m \pi(z) = \<z, \nabla \pi(z)\>$ for all $z = (x,y)\in \R^2$. 
Assuming without loss of generality that $\alpha \neq 0$,
we therefore obtain
$$
m \pi(x,y) =  l^s (x\pi_1 +y \pi_2) = l^s\left( (\alpha x+\beta y)\frac {\pi_1} \alpha + \frac y \alpha (\alpha \pi_2 - \beta \pi_1)\right)
$$
which shows that $\pi$ is a multiple of $l^{s+1}$, hence of $l^{t_m}$.  

Conversely if $\pi$ is a multiple of $l^{t_m}$ then $\partial_x \pi$ and $\partial_y \pi$ are both multiples of $l^{t_m-1}$. Since $2(t_m-1)\geq m$ the polynomial $|\nabla \pi|^2$ is a multiple of $l^m$ which concludes the proof.
\sq

\section{Proof of the main result}
\label{secProofs}
This section is devoted to the proof of Theorem \ref{optiTheorem}. We thus consider a fixed bounded polygonal domain $\Omega\subset \R^2$, an integer $m\geq 2$, an exponent $1\leq p < \infty$ and a function $f\in C^m(\overline \Omega)$.

The Taylor development of $f$ close to a point $z\in \Omega$ is given by two polynomials $\mu_z\in \P_{m-1}$ and $\pi_z\in \H_m$: $f(z+h) = \mu_z(h)+\pi_z(h)+o(|h|^m)$, where $h\in \R^2$ is small. The Taylor development of the function $\nabla f : \Omega\to \R^2$ close to a point $z\in \Omega$ is obtained by derivation of the previous one:
$$
\nabla f(z+h) = \nabla \mu_z(h)+ \nabla \pi_z(h)+ o(|h|^{m-1}),
$$
and the corresponding Taylor formula reads as follows: for any $z,h$ such that $[z,z+h]\subset \Omega$
\be
\label{TaylorInt}
\nabla f(z+h) = \nabla \mu_z(h)+ (m-1)\int_{t=0}^1 \nabla \pi_{z+th} (h) (1-t)^{m-2} dt.
\ee
We denote by $\omega$ the modulus of continuity of the function $z\in \Omega\mapsto \nabla \pi_z\in \H_{m-1}^2$: for each $r>0$
\be
\label{defomega}
\omega(r) := \sup \left\{ \|\nabla \pi_z - \nabla \pi_{z'}\| \ssep z,z'\in \Omega \text{ s.t. } |z-z'|\leq r\right\},
\ee
where we recall that $\|\mu\| = \sup\{ |\mu(z)|/|z|^{m-1}\ssep z\neq 0\}$ for all $\mu\in \H_{m-1}^2$, see \iref{defNormH2}. 

Our first lemma is an estimation on a single triangle of the gradient interpolation error of $f$. %, that will be used both in this section and the next.
\begin{lemma}
\label{lemmaDiff}
There exists a constant $\CDiff=\CDiff(f,m)$ such that the following holds.
\begin{enumerate}[i.]
\item
For any triangle $T\subset \overline \Omega$, and any $z\in T$ one has
\be
\label{errorDiff}
\left | \|\nabla (f-\interp_T^{m-1} f)\|_{L^p(T)} - \|\nabla (\pi_z - \interp_T^{m-1} \pi_z)\|_{L^p(T)} \right| \leq \CDiff \omega(\diam(T)) |T|^\frac 1 \tau \rho(T)^{\frac {m+1} 2} 
\ee
\item
For any triangle $T\subset \overline \Omega$ one has 
\be
\label{estimBd}
\|\nabla (f-\interp_T^{m-1} f)\|_{L^p(T)} \leq C_\star |T|^\frac 1 p \diam(T)^{m-1} S(T).
\ee
\end{enumerate}
\end{lemma}

\begin{proof}
The point $z$ is fixed throughout this proof. We define a function $g\in C^m(\overline \Omega)$ by the equality
$$
g(z+h) := f(z+h) - \mu_z(h)- \pi_z(h),
$$
for all $h$ such that $z+h\in \overline \Omega$. We have 
\begin{eqnarray*}
\left | \|\nabla (f-\interp_T^{m-1} f)\|_{L^p(T)} - \|\nabla (\pi_z - \interp_T^{m-1} \pi_z)\|_{L^p(T)} \right| 
&\leq& \|\nabla ((f-\pi_z) -\interp_T^{m-1} (f-\pi_z))\|_{L^p(T)}\\
&=& \|\nabla (g -\interp_T^{m-1} g)\|_{L^p(T)}\\
& \leq & |T|^\frac 1 p  \|\nabla (g -\interp_T^{m-1} g)\|_{L^\infty(T)}\\
& \leq & C \rho(T) |T|^\frac 1 p \|\nabla g\|_{L^\infty(T)},
\end{eqnarray*}
where we used the reverse triangle inequality in the first line, the translation invariance \iref{transInv} and the equality  $\mu_z = \interp_T^{m-1}\mu_z$ in the second line, Jensen's inequality in the third line, and Lemma \ref{lemmaS} in the last line combined with the inequality $1\leq S(T) \leq \rho(T)$, see \iref{ineqRhoS}.
We have
$$
\nabla g(z+h) = \nabla f(z+h) - \nabla \mu_z(h) - \nabla \pi_z(h) = (m-1)\int_{t=0}^1 (\nabla\pi_{z+th}(h)-\nabla\pi_z(h))(1-t)^{m-2}dt,
$$
therefore 
$
\|\nabla g(z+h)\| \leq |h|^{m-1} \omega(h). 
$
Furthermore if $z\in T$ and $z+h\in T$, then $|h|^2 \leq \diam(T)^2 \leq (4/|\TEq|) |T| \rho(T)$ according to Lemma \ref{lemmaDiamRho}. This concludes the proof of Point i, provided that $C_\star \geq C(4/|\TEq|)^\frac{m-1} 2$.

We now turn to the proof of Point ii, and for that purpose we consider a fixed point $z\in T$. Changing our previous notation we denote by $g\in C^m(\overline \Omega)$ the function defined by $g(z+h) := f(z+h) - \mu_z(h)$ for all $h$ such that $z+h\in \overline \Omega$.
We obtain 
$$
\|\nabla (f-\interp_T^{m-1} f)\|_{L^p(T)} = \|\nabla (g-\interp_T^{m-1} g)\|_{L^p(T)} \leq |T|^\frac 1 p\|\nabla (g-\interp_T^{m-1} g)\|_{L^\infty(T)} \leq C|T|^\frac 1 p \|\nabla g \|_{L^\infty(T)},
$$
where we used successively that the interpolation operator reproduces the elements of $\P_{m-1}$, Jensen's inequality, and Lemma \ref{lemmaS}. On the other hand we obtain using \iref{TaylorInt}.
$$
\nabla g(z+h) = \nabla f(z+h) -\nabla \mu_z(h) = (m-1)\int_{t=0}^1 \nabla \pi_{z+th}(h) (1-t)^{m-2} dt,
$$
hence $|\nabla g(z+h)| \leq |h|^{m-1} \sup \{\|\nabla \pi_z\| \ssep z\in \Omega\}$. If $z\in T$ and $z+h\in T$, then $|h| \leq \diam(T)$, which concludes the proof of \iref{estimBd} provided that $C_\star \geq  C\sup \{\|\nabla \pi_z\| \ssep z\in \Omega\}$.
\end{proof}

\subsection{The lower error estimate \iref{lowerEstim}}
\label{secLowerEstim}

Under the hypotheses of Lemma \ref{lemmaDiff}, and recalling from \iref{transInv} that $\|\nabla (\pi- \interp_T^{m-1} \pi)\|_{L^p(T)} \geq |T|^{\frac 1 \tau} L_{m,p}(\pi)$, we obtain
\be
\label{estimLowerLoc}
\|\nabla (f-\interp_T^{m-1} f)\|_{L^p(T)} \geq |T|^\frac 1 \tau \left(L_{m,p}(\pi_z) - \CDiff\omega(\diam T) \rho(T)^{\frac{m+1} 2}\right).
\ee

We consider an admissible sequence of triangulations $(\cT_N)_{N\geq N_0}$. For all $N\geq N_0$, all
$T\in \cT_N$ and all $z\in T$, we define $\phi_N(z) := |T|$ and 
$$
\psi_N(z) := \left(L_{m,p}(\pi_z) - \CDiff\omega(\diam(T)) \rho(T)^{\frac{m+1} 2}\right)_+,
$$
where $\lambda_+ := \max\{\lambda,0\}$.
Holder's inequality $\int f_1 f_2 \leq \|f_1\|_{p_1} \|f_2\|_{p_2}$, applied to the functions
$
f_1 = \phi_N^{\frac {(m-1) \tau} 2} \psi_N^\tau$  and  $f_2 = \phi_N^{-\frac {(m-1) \tau} 2},
$
and the exponents
$
p_1 = \frac p \tau \ \text{ and } \ p_2 = \frac 2 {(m-1)\tau},
$ 
yields 
\be
\label{holderPsi}
\int_\Omega \psi_N^\tau \leq \left(\int_\Omega \phi_N^{\frac{(m-1)p} 2} \psi_N^{p}\right)^{\frac \tau p} \left(\int_\Omega \phi_N^{-1}\right)^{\frac {(m-1)\tau} 2}.
\ee
Furthermore for any $T\in \cT_N$ and any $z\in T$, we obtain using \iref{estimLowerLoc} 
$$
 \phi_N(z)^{\frac{(m-1)p} 2} \psi_N(z)^p
=|T|^{\frac p \tau-1}  \psi_N(z)^p \leq\frac 1 {|T|} \|\nabla(f-\interp_T^{m-1} f)\|_{L^p(T)}^p,
$$
hence 
$$
\int_\Omega \phi_N^{\frac{(m-1)p} 2} \psi_N^{p}\leq \sum_{T\in \cT_N} \frac 1 {|T|}  \int_T  \|\nabla(f-\interp_T^{m-1} f)\|_{L^p(T)}^p= \|\nabla(f-\interp_{\cT_N}^{m-1} f)\|_{L^p(\Omega)}^p.
$$
Elevating \iref{holderPsi} to the power $\frac 1 \tau$, injecting the above estimate and observing that $\int_\Omega \phi_N^{-1} = \# \cT_N\leq N$, 
 we thus obtain
\be
\label{upperPsi}
 \|\psi_N\|_{L^\tau(\Omega)} \leq \|\nabla(f-\interp_{\cT_N}^{m-1} f)\|_{L^p(\Omega)} N^{\frac {m-1} 2}.
\ee

\noindent
Since the sequence $\seqT$ is admissible, there exists a constant $C_A>0$ such that for all $N$ and all $T\in \cT_N$ we have $\diam(T)\leq C_AN^{-\frac 1 2}$. 
We introduce a subset of $\cT'_N\subset \cT_N$ which gathers the most degenerate triangles
$$
\cT'_N = \{ T\in \cT_N \sep \rho(T)\geq \omega(C_AN^{-\frac 1 2})^{\frac{-1}{m+1}}\},
$$
where $\omega$ is defined by \iref{defomega}. We denote by $\Omega'_N$ the portion of $\Omega$ covered by $\cT'_N$.
For all $z\in \Omega\sm\Omega'_N$ one has %, recalling from \iref{ineqRhoS} that $\rho \geq S$, 
%we obtain 
$$
\psi_N(z)\geq L_{m,p}(\pi_z) -\CDiff\sqrt{\omega(C_A N^{-\frac 1 2})}.
$$
Hence, with the convention $L_{m,p}(d^m f(z)/m!) := L_{m,p}(\pi_z)$,
\begin{eqnarray*}
 \|\psi_N\|_{L^\tau(\Omega)}^\tau & \geq& \left \|\left(L_{m,p}\left(\frac{d^m f}{m!} \right) -\CDiff\sqrt{\omega(C_A N^{-\frac 1 2})}\right)_+\right\|_{L^\tau(\Omega\sm\Omega'_N)}^\tau\\
 & \geq &
 \left \|\left(L_{m,p}\left(\frac{d^m f}{m!} \right) -\CDiff\sqrt{\omega(C_A N^{-\frac 1 2})}\right)_+\right\|_{L^\tau(\Omega)}^\tau
 -C^\tau |\Omega'_N|,
 \end{eqnarray*}
where $C:=\max \{L_{m,p}(\pi_z)\ssep z\in \Omega\}<\infty$. We next observe
that $|\Omega'_N|\to 0$ as $N\to +\infty$: indeed
for all $T\in \cT'_N$ we obtain using Lemma \ref{lemmaDiamRho}
$$
|T| \leq |\TEq|\diam(T)^2 \rho(T)^{-1} \leq |\TEq|C_A ^2 N^{-1} \omega(C_A N^{-\frac 1 2})^{\frac 1 {m+1}}.
$$
Since $\#(\cT'_N)\leq N$, we obtain $|\Omega'_N|\leq |\TEq|C_A^2 \omega(C_A N^{-\frac 1 2})^{\frac 1 {m+1}}$, which tends to $0$ as $N\to \infty$. Therefore %We thus obtain
$$
\liminf_{N\to \infty} \|\psi_N\|_{L^\tau(\Omega)}
\geq \lim_{N\to \infty}   \left \|\left(L_{m,p}\left(\frac{d^m f}{m!} \right) -\sqrt{\omega(C_A N^{-\frac 1 2})}\right)_+\right\|_{L^\tau(\Omega)}
 = \left\|L_{m,p}\left(\frac{d^m f}{m!} \right)\right\|_{L^\tau(\Omega)}.
$$
Combining this result with \iref{upperPsi} we conclude the proof of the announced estimate \iref{lowerEstim}.

\subsection{A triangulation containing small periodic patches}

%\section{The upper estimate}
%\label{secUpperEstim}

%This section is devoted to the proof of the upper estimate \iref{upperEstimEps} of Theorem \ref{optiTheorem}. We first describe the construction of a sequence of triangulations, which consists of the aggregation of small periodic patches. We then use this construction to obtain the announced estimates.

%The polygonal bounded domain $\Omega$, the integer $m$, the exponent $p$ and the function $f\in C^m(\overline \Omega)$ are fixed throughout this section.
%*)
This subsection describes the construction of some triangulations by the aggregation of small periodic patches, which is a preliminary step for the proof of the estimate \iref{upperEstimEps} of Theorem \ref{optiTheorem} in the next subsection.
The design of these triangulations is related to the anterior works \cite{BBLS,BLM}, yet it is adapted in order to keep under control the measure of sliverness of the elements. 
The triangulations considered in this subsection are denoted by the letter $\cP$, instead of $\cT$, in order to avoid conflicts of notations with the next subsection.

Our first lemma describes a family of meshes of a triangle $R$, which consist for the largest part of a periodic tiling based on a triangle $T$ scaled down by a factor $1/n$, except for a few elements close to $\partial R$. 

\begin{lemma}
\label{lemmaTile}
Let $R$ and $T$ be two triangles. There exists a family $(\cP_{T,R,n})_{n\geq 1}$, 
of conforming triangulations of $R$, and a constant $C_{R,T}$, such that the following holds.
\begin{enumerate}
\item (Tights bounds on the cardinality of $\cP_{T,R,n}$ and the diameter of its elements)
\be
\label{cardDiamTile}
\lim_{n\to \infty} \frac {\# (\cP_{T,R,n})} {n^2} = \frac {|R|}{|T|}
% \stext{ and } \sup_{n\geq 1}  \; \( n \max_{T'\in \cP_{T,R,n}} \diam(T') \)\leq 3 (\diam(T)+  \diam(R)).
\stext{ and }  \max_{T'\in \cP_{T,R,n}} \diam(T') \leq \frac 3 n (\diam(T)+  \diam(R)). % \text{ for all } n \geq 1.
\ee

\item (Conformity) The vertices of $\cP_{T,R,n}$ on the boundary of $R$ are exactly those of the form $\frac k n a+ (1-\frac k n) b$, where $0\leq k \leq n$ and $a,b$ are vertices of $R$.

\item (Control of the boundary elements) Denote by $R_{T,n}\subset R$ the union of all the elements of $\cP_{T,R,n}$ which are \emph{not} of the following form: the image of $T$ by a map of the form $z\mapsto z_0+\sigma z/n$, for some $z_0\in \R^2$, $\sigma\in \{-1,1\}$. Then for all $n \geq 1$
\be
\label{boundT}
|R_{T,n}| \leq \frac {C_{R,T}} n \stext{ and } \max_{T'\in \cP_{R,T,n}} S(T') \leq C_{R,T}. 
\ee
\end{enumerate}
\end{lemma}

\begin{proof}
See Appendix. 
\end{proof}

%We introduce the concept of \emph{local shape specification}, which encodes a  and in that sense is loosely related to the concept of Riemannian metrics discussed in \S\ref{subsecConst}.

%We introduce the notion of \emph{local shape specification}, which is the counterpart in this section of the plays a r
We next introduce the concept of \emph{local shape specification}, on the domain $\Omega$.
\begin{definition}
\label{defLocShape}
A \emph{local shape specification} is a (possibly discontinuous) map $y\in \overline \Omega \mapsto T_y$ which associates a triangle $T_y$ to each point $y$ in the closure of $\Omega$, and which satisfies the following properties.
\begin{itemize}
\item The volume map $y\in \overline \Omega \mapsto |T_y|\in \R_+^*$ is \emph{continuous and positive}. % map.
\item The measure of degeneracy $y \in \overline \Omega \mapsto \rho(T_y) \in \R_+$ is \emph{uniformly bounded}.
\end{itemize}
\end{definition}

The next proposition describes a sequence of triangulations adapted in a certain sense to a given local shape specification. This can be compared to the construction, evoked in \S\ref{subsecConst}, of a triangulation adapted to a given Riemmannian metric in the sense of \iref{adaptTH}.
Given a point $y\in \R^2$ and compact set $T$ we define $d_\cH(y, T) := \max \{|y-z| \ssep z\in T\}$, which is the Haussdorf distance separating the sets $\{y\}$ and $T$. 
\begin{prop}
\label{propTile}
Let $y\mapsto T_y$ be a local shape specification. There exists a sequence $(\cP_n)_{n \geq 2}$ of triangulations of $\Omega$, a sequence $(\delta_n)_{n \geq 2}$ of positive reals converging to $0$, and a constant $C_a$, satisfying the following properties.
\begin{enumerate}
\item (Tight bounds on the cardinality of $\cP_n$ and the diameter of its elements)
\be
\label{limCardPn}
\lim_{n \to \infty} \frac{\#(\cP_n)} {n^2} = \int_\Omega \frac {dy} {|T_y|} \stext{ and } 
\max_{T\in \cP_n} \diam(T) \leq \frac{C_a} n.
%\sup_{n \geq N_0} \left(n \max_{T \in \cP_n} \diam(T)\right) < \infty.
\ee
\item (Control of the boundary elements)
Denote by $\Omega_n\subset \Omega$ the union of all the elements $T\in \cP_n$ which are \emph{not} of the following form: the image of $T_y$, for some $y\in \Omega$ such that $d_\cH(y,T) \leq \delta_n$, by a map of the form $z\mapsto z_0+ \sigma z/n$, for some $z_0\in \R^2$, $\sigma\in \{-1,1\}$. Then for all $n \geq 2$
$$
|\Omega_n| \leq \frac{C_a \ln n} n \stext{ and } \max_{T\in \cP_n} S(T) \leq  C_a \ln n.
$$
\end{enumerate}
\end{prop}

\begin{proof}
We consider a triangulation $\cR^1$ of the polygonal domain $\Omega$ of minimal cardinality $N_0$. For each $k \geq 1$ we denote by $\cR^k$ the triangulation of $\Omega$ of cardinality $k^2 N_0$ obtained by uniformly subdividing the elements of $\cR^1$ into $k^2$ sub triangles.

For each $n\geq 1$ we denote by $\cR^k_n$ the triangulation of $\Omega$ obtained as the union of the triangulations $\cP_{R,T,n}$ described in Lemma \ref{lemmaTile}, where $R\in \cR^k$ and $T=T_{z_R}$ is the triangle specified by the local shape specification at the barycenter of $R$.
Point 2 of Lemma \ref{lemmaTile} guarantees that this triangulation is \emph{conforming}: there is no hanging node at the interfaces of the triangles $R\in \cR^k$. Our next observation is that for each $k,n \geq 1$, using Lemma \ref{lemmaDiamRho}
$$
n \max_{T\in \cR_n^k} \diam(T) \leq 3 \max_{R \in \cR^k} (\diam(R)+ \diam(T_{z_R})) \leq  3 \diam(\Omega)+\sup_{y\in \Omega} \sqrt{(4/|\TEq|)|T_y| \rho(T_y)},
$$
which is finite and independent of $k$ and $n$. For each $k,n\geq 1$ we define a real $\delta(n,k)$ by the equality
$$
\frac{\#(\cR_n^k)}{n^2} = \int_\Omega \frac{dy}{|T_y|} + \delta(n,k).
$$
Using \iref{cardDiamTile} we obtain 
$$
\lim_{n\to \infty} \delta(n,k) = \delta(k) := \sum_{R\in \cR^k} \frac {|R|}{|T_{z_R}|} -  \int_\Omega \frac{dy}{|T_y|}.
$$
Note that $\delta(k) \to 0$ as $k \to \infty$ since the map $y \mapsto T_y$ is continuous.

We denote by $\Omega_n^k$ the union of the sets $R_{T,n}$ described in Lemma \ref{lemmaTile}, for $R\in \cR^k$ and $T=T_{z_R}$, and by $C_k$ the sum of the corresponding constants $C_{R,T}$. %, for the same collection of pairs $(R,T)$, .
We obtain using \iref{boundT} that for all $n \geq 1$
$$
|\Omega_n^k|\leq \frac{C_k} n \stext{ and } \max_{T\in \cR_n^k} S(T) \leq C_k.
$$
We finally choose a sequence $(k(n))_{n \geq 2}$, such that $k(n)\to \infty$ as $n \to \infty$, and which increases ``slowly'' in the following sense: we require that for $n$ sufficiently large one has $\delta(n,k(n)) \leq \delta(k(n))+\frac 1 n$ and $C_{k(n)} \leq \ln n$.
Defining $\cP_n := \cR_n^{k(n)}$ and $\delta_n := \diam(\Omega)/k(n)$ we obtain the announced result.
\end{proof}

\subsection{The upper error estimate \iref{upperEstimEps}}

Throughout this section we consider a fixed real $M \geq 1$, and we introduce the collection of triangles 
$$
\bT_M := \{ T \ssep |T|=1,\, \rho(T)\leq M, \, z_T=0\}
$$
which is compact for the Haussdorf distance.
We introduce a variant $L_M$ of the shape function $L_{m,p}$ defined as follows: for each $\pi \in \H_m$
\be
\label{defLM}
L_M(\pi) := \min_{T\in \bT_M} \|\nabla (\pi-\interp_T^{m-1} \pi)\|_{L_p(T)}.
\ee
For any fixed $\pi \in \H_m$ the map $T\mapsto \|\nabla(\pi - \interp_T^{m-1} \pi)\|_{L^p(T)}$
is continuous with respect to the Haussdorff distance on the set of all triangles. Hence there exists a minimizing triangle, that we denote by $T(\pi)\in \bT_M$, for the optimization problem appearing in \iref{defLM}.

Since all norms are equivalent on the finite dimensional space $\H_m$, there exists a constant $C_M$ such that for all $\pi\in \H_m$ 
$$
\sup_{T \in \bT_M} \|\nabla (\pi - \interp_T^{m-1} \pi)\|_{L^p(T)} \leq C_M \|\nabla \pi\|.
$$
The function $L_M$ is defined as the infimum of a family of $C_M$-Lipschitz functions on $\H_m$, hence is also $C_M$-Lipschitz: $|L_M(\pi) - L_M(\pi')|\leq C_M \|\nabla \pi-\nabla \pi'\|$. Furthermore for each $\pi, \pi'\in \H_m$, one has since $T(\pi) \in \bT_M$
\be
\label{changeT}
\| \nabla (\pi'-\interp_{T(\pi)}^{m-1} \pi')\|_{L^p(T(\pi))} \leq \| \nabla (\pi-\interp_{T(\pi)}^{m-1} \pi)\|_{L^p(T(\pi))} + C_M\|\pi-\pi'\| = L_M(\pi) +  C_M\|\pi-\pi'\|.
\ee

We consider the following \emph{local shape specification}, see Definition \ref{defLocShape},
$$
y \mapsto T_y := (L_M(\pi_y)+ M^{-1})^{-\frac \tau 2} T(\pi_y).
$$
In other words $T_y$ is the isotropic scaling of the triangle $T(\pi_y)$ by the factor $(L_M(\pi_y)+ M^{-1})^{-\frac \tau 2}$. We thus have $|T_y| = (L_M(\pi_y)+ M^{-1})^{-\tau}$, which depends continuously on $y\in \overline \Omega$, and $\rho(T_y) = \rho(T(\pi_y))\leq M$.

Proposition \ref{propTile} describes a sequence $(\cP_n)_{n \geq 2}$ of triangulations of $\Omega$ attached to this local shape specification, as well as a sequence $(\delta_n)_{n \geq 2}$ of positive numbers, a sequence $(\Omega_n)_{n \geq 2}$ of subdomains of $\Omega$, and a constant $C_a$.
We recall that for all $n \geq 2$ and all $T\in \cP_n$
$$
\diam(T) \leq C_a/n, \quad S(T) \leq C_a \ln n, \quad  |\Omega_n| \leq C_a (\ln n)/n.
$$

Using Point ii of Lemma \ref{lemmaDiff} we obtain for each $n \geq 1$ and each triangle $T \in \cP_n$
$$
\|\nabla (f-\interp_T^{m-1} f)\|_{L^p(T)} \leq C_\star C_a^m |T|^\frac 1 p (\ln n)/n^{m-1}.
$$
Summing up the $p$-th power of the contributions of all the triangles $T\in \cP_n$ such that $T\subset \Omega_n$, we obtain 
$$
\|\nabla (f-\interp_T^{m-1} f)\|_{L^p(\Omega_n)} \leq C_\star C_a^m |\Omega_n|^\frac 1 p \frac{\ln n}{n^{m-1}}\leq  \frac{C_\star C_a^{m+1}(\ln n)^{1+\frac 1 p}n^{-\frac 1 p}}{n^{m-1}} =: \frac{\ve_n}{n^{m-1}},
$$
and we observe that $\ve_n\to 0$ as $n \to \infty$.

We now turn to the contribution of $\Omega\sm \Omega_n$ to the error, and for that purpose we consider a triangle $T$ which is the image of the triangle $T_y$, for some $y\in \Omega$ such that $d_\cH(y,T) \leq \delta_n$, by a map of the form $z\mapsto z_0+ \sigma z/n$, for some $z_0\in \R^2$, $\sigma\in \{-1,1\}$. We have for any $z\in T$, using Point i of Lemma \ref{lemmaDiff}
\begin{eqnarray*}
\|\nabla (f-\interp_T^{m-1} f)\|_{L^p(T)} &\leq& \|\nabla (\pi_z-\interp_T^{m-1} \pi_z)\|_{L^p(T)} + C_\star \omega(C_a/n) |T|^\frac 1 \tau M^{\frac {m+1} 2}\\
& = & |T|^{\frac 1 \tau} \left(\|\nabla (\pi_z - \interp_{T(\pi_y)}^{m-1} \pi_z)\|_{L^p(T(\pi_y))} +C_\star \omega(C_a/n) M^{\frac {m+1} 2}\right)\\
& \leq & \frac{L_M(\pi_y)+ C_M \omega(\delta_n) + C_\star M^\frac{m+1} 2\omega(C_a/n)}{ n^{\frac 2 \tau}(L_M(\pi_y)+ M^{-1})},
\end{eqnarray*}
where we used %Lemma \ref{lemmaDiff}, Point i., in the first line, 
the invariance property \iref{transInv} in the second line, and \iref{changeT} in the last line.
It follows that there exists an integer $n_0$ such that  $\|\nabla (f-\interp_T^{m-1} f)\|_{L^p(T)}\leq n^{-\frac 2 \tau}$ for all $n \geq n_0$ and all $T\in \cP_n$.
Hence 
\begin{eqnarray*}
\|\nabla (f-\interp_{\cP_n}^{m-1} f)\|_{L^p(\Omega)} &\leq& \|\nabla (f-\interp_{\cP_n}^{m-1} f)\|_{L^p(\Omega\sm\Omega_n)} +  \|\nabla (f-\interp_{\cP_n}^{m-1} f)\|_{L^p(\Omega_n)} \\
& \leq & \#(\cP_n)^\frac 1 p n^{-\frac 2 \tau} + \ve_n n^{-(m-1)},
\end{eqnarray*}
and therefore using \iref{limCardPn}
\begin{IEEEeqnarray*}{rCl}
\limsup_{n\to \infty} \#(\cP_n)^{\frac{m-1} 2} \|\nabla(f-\interp_{\cP_n}^{m-1} f)\|_{L^p(\Omega)} &\leq& \lim_{n \to \infty} \frac{\#(\cP_n)^\frac 1 \tau}{n^\frac 2 \tau} + \frac{\#(\cP_n)^{\frac {m-1}2} \ve_n}{n^{m-1}}\\
 &=& \left(\int_\Omega \frac {dy}{|T_y|}\right)^\frac 1 \tau+0\\
 &=& \left(\int_\Omega (L_M(\pi_y)+M^{-1})^\tau dy\right)^\frac 1 \tau.
\yesnumber
\label{limErrPn}
\end{IEEEeqnarray*}

Observe that $L_M(\pi_y)$ converges \emph{decreasingly} for each $y\in \Omega$ to $L_{m,p}(d^m f(y)/m!) := L_{m,p}(\pi_y)$ as $M\to \infty$. 
Given any fixed $\ve>0$, and using standard results on the convergence of integrals, we may therefore choose $M = M(\ve)$ sufficiently large such that 
\be
\label{LMLmp}
\left(\int_\Omega (L_M(\pi_y)+M^{-1})^\tau dy\right)^\frac 1 \tau \leq \left\| L_{m,p}\left(\frac {d^m f}{m!}\right)\right\|_{L^\tau(\Omega)}+\ve.
\ee
We denote by $N_0$ the minimal cardinality of a triangulation of the polygonal domain $\Omega$, and we assume without loss of generality that $\#(\cP_2) = N_0$.
For each $N \geq N_0$, we denote by $n(N)$ the largest integer such that $\#(\cP_{n(N)}) \leq N$, and we set $\cT_N^\ve := \cP_{n(N)}$. Observing that $\#(\cT_N^\ve)/N \to 1$ as $N \to \infty$, and combining \iref{limErrPn} with \iref{LMLmp}, we obtain the announced result \iref{upperEstimEps}:
$$
\limsup_{N \to \infty} N^{\frac {m-1} 2} \|\nabla (f-\interp_{\cT_N^\ve}^{m-1}f)\|_{L^p(\Omega)} = \limsup_{n\to \infty} \#(\cP_{n(N)})^{\frac{m-1} 2} \|\nabla(f-\interp_{\cP_{n(N)}}^{m-1} f)\|_{L^p(\Omega)} \leq  \left\| L_{m,p}\left(\frac {d^m f}{m!}\right)\right\|_{L^\tau(\Omega)}\hspace{-4mm}+\ve.
$$
The admissibility of the sequence $(\cT_N^\ve)_{N \geq N_0}$ of triangulations immediately follows from \iref{limCardPn}.
\section*{Conclusion}

In this paper, we have introduced asymptotic estimates for the
finite element interpolation error measured in the $W^{1,p}$ semi-norm, when the mesh is 
optimally adapted to a function of two variables and the degree of interpolation $m-1$ is arbitrary.
The approach used is an adaptation of the ideas developped in \cite{Mi} for the $L^p$ interpolation error, and leads to asymptotically sharp error estimates, exposed in Theorems \ref{mainTheorem} and \ref{optiTheorem}.
These estimates involve a shape function $L_{m,p}$ which generalises the determinant which appears in estimates for piecewise linear interpolation. The shape function has equivalents of polynomial 
form for all values of $m$, as established in \S\ref{secPolEq}. 
Up to a fixed multiplicative constant, our estimates can therefore be written 
under analytic form in terms of the derivatives of the function to be approximated. 

Future efforts will be devoted to the extension of these results to functions defined on a domain of dimension $d>2$, which is partially done in Chapter 3 of \cite{thesisJM}. Another challenge left open is the development of an anisotropic mesh generator with guarantees on the maximal angle of the elements, as evoked in \S \ref{secPractMesh}, which could allow to apply the results of this paper in the context of adaptive mesh refinement for numerical simulations.

\appendix
\begin{center}
    {\bf APPENDIX}
  \end{center}

\section{Proof of Lemma \ref{lemmaTile}}

We construct the triangulation $\cP_{T,R,n}$ of $R$ as the union of three components: $\cP_{T,R,n} = \cB_n \cup \cL_n\cup \cI_n$. The elements of $\cB_n$ cover the \emph{boundary} of $R$, while the elements of $\cI_n$ cover most of its \emph{interior} and are included in $\overline {R\sm R_{T,n}}$. The elements of $\cL_n$ play the role of a \emph{layer} between $\cI_n$ and $\cB_n$. Throughout this proof we denote by $C=C(R,T)$ a generic constant independent of $n$, which may change from one occurrence to the next.

We introduce the homothetic contraction $R_n$ of the triangle $R$ by the factor $1-n^{-1}$, with the same barycenter (i.e.\ the image of $R$ by the map $z\mapsto z_R+ (1-n^{-1})(z-z_R)$). One easily checks that for all $z\in \partial R$
\be
\label{distToRn}
%n \inf \{ |z-y| \ssep y\in R_n\} \leq \frac 2 3\diam(T). 
\frac 2 3 \frac {|R|}{\diam(R)} \leq n \, d(z,R_n) \leq \frac 2 3\diam(R), 
\ee
where $d(z,E) := \inf\{|z-e| \ssep e\in E\}$ for any $z\in \R^2$ and any $E\subset \R^2$.

\begin{figure}
	\centering
		\includegraphics[width=4cm,height=4cm]{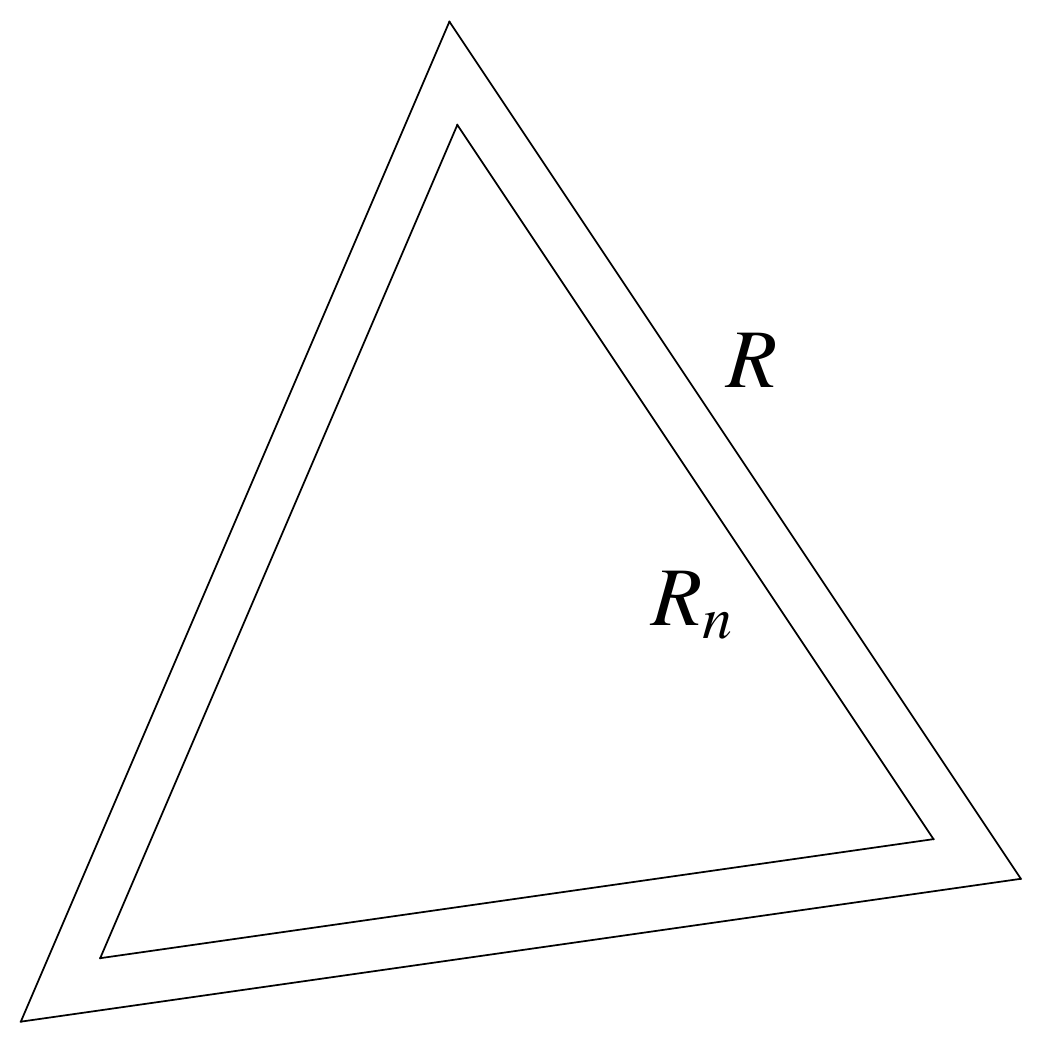} %{Illustrations/RandRn.eps}
%		\hspace{1cm}
		\includegraphics[width=4cm,height=4cm]{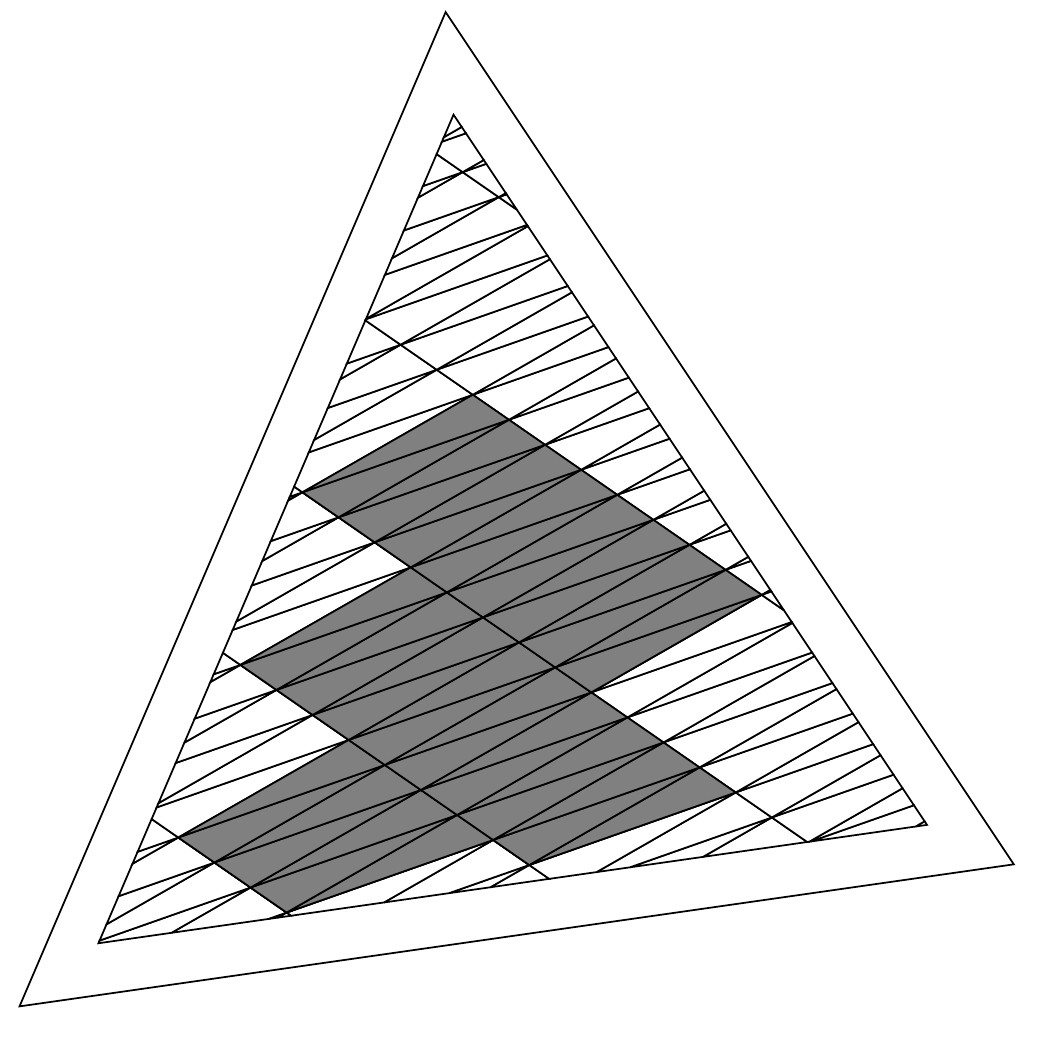} %{Illustrations/TilingPure.eps}
		\includegraphics[width=3.5cm,height=3.5cm]{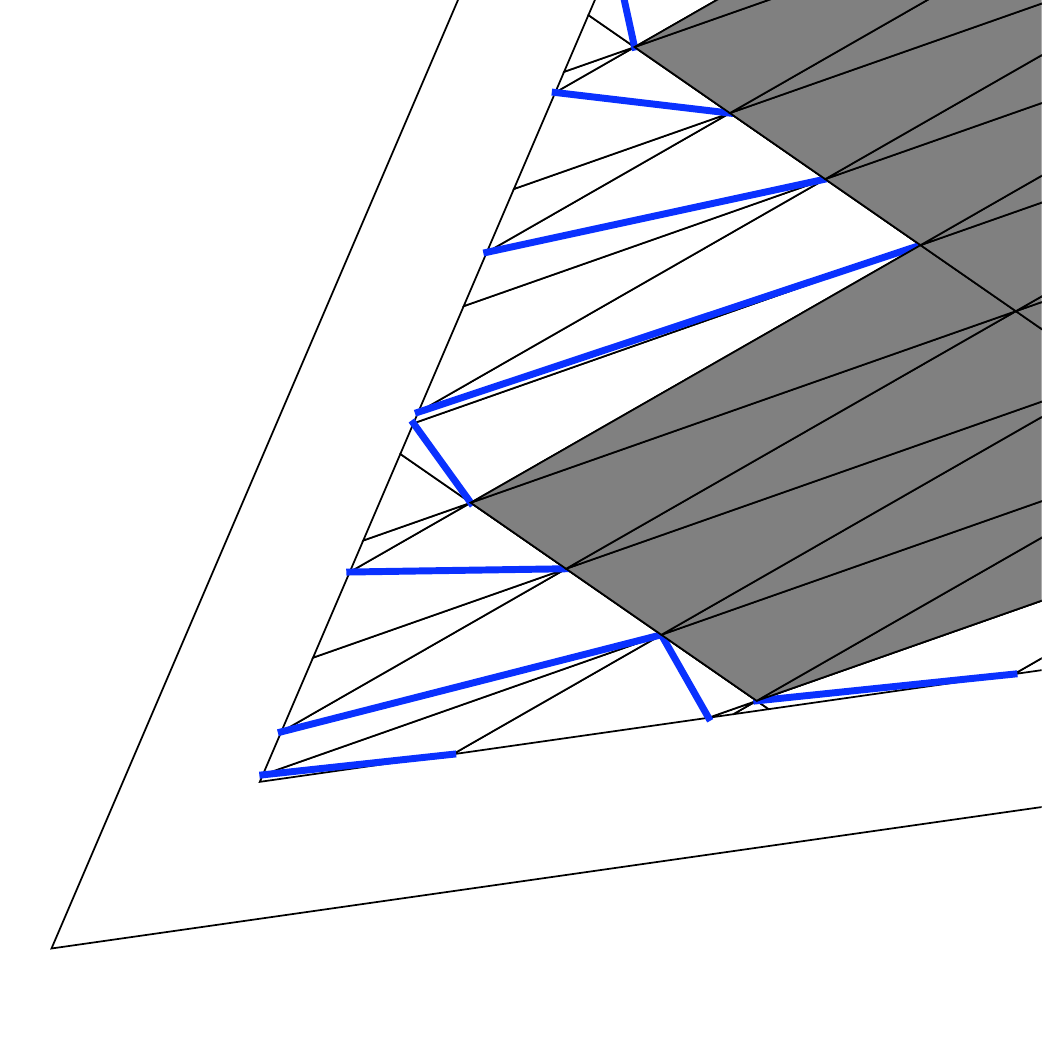} %{Illustrations/TilingBlue.eps}
%		\hspace{1cm}
		\includegraphics[width=3.5cm,height=3.5cm]{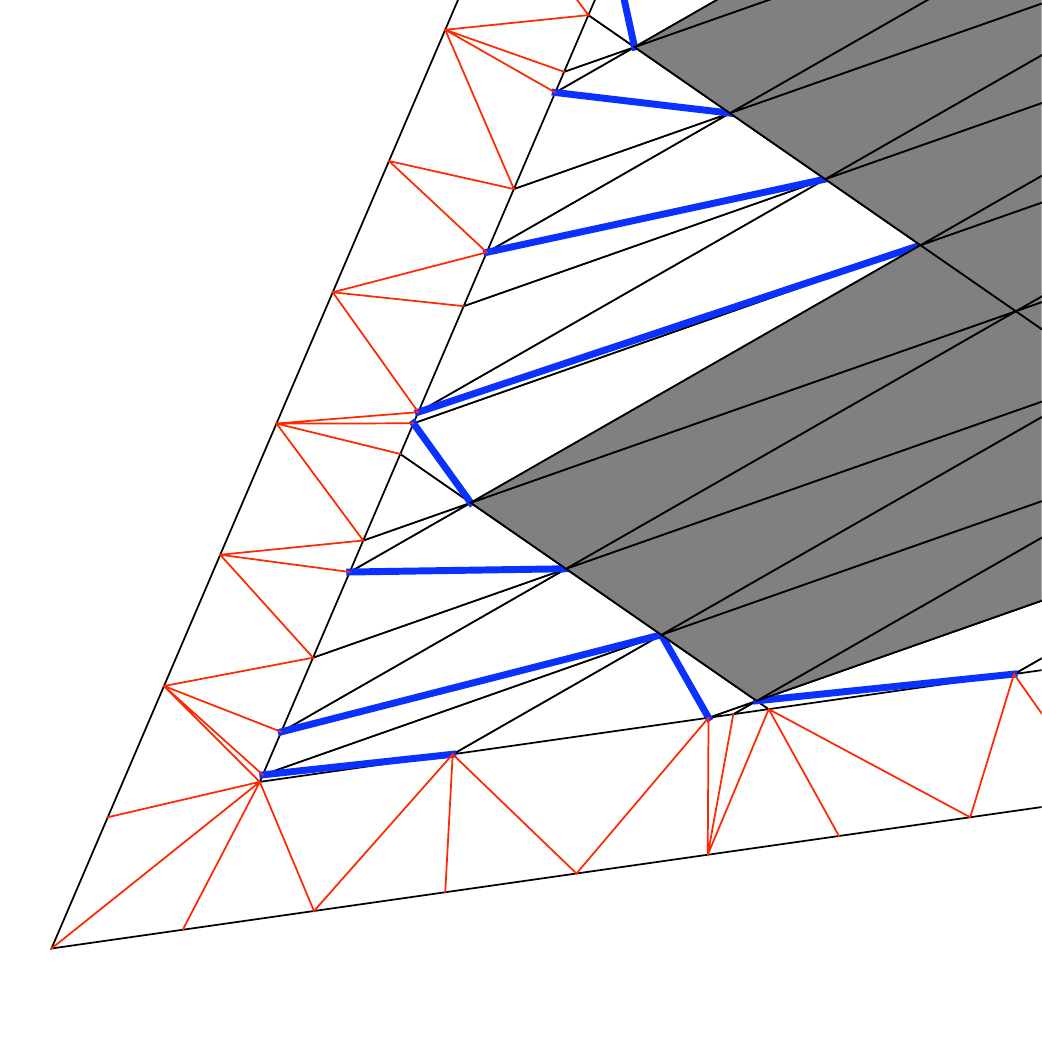} %{Illustrations/TilingRed.eps}

	\caption{The triangles $R$ and $R_n$ (left). The collections $\cI_n$ (gray) of triangles and $\cL_n^0$ (white) of convex polygons (center left). The collection $\cL_n$ of triangles is obtained by triangulating the elements of $\cL_n^0$ (center right). The collection $\cB_n$ of triangles covers $R\sm R_n$ (right).}
%	The partition $\cP'_n$ of $R_n$ (center left). Detail of the partition $\cP''_n$ of $R_n$ (center right). The partition $\cP_n = \cP_{T,n}(R) = \cP''_n \cup \widetilde \cP_n$ of $R$ (right).} 
	\label{figAppendix}
\end{figure}

Let $z_0$ be a vertex of the triangle $T$, and let $u,v$ be two of its edge vectors. 
For each $n \geq 1$ we denote by $\cI_n^0$ the periodic tiling of $\R^2$ built of the images of $T$ by the maps $z\mapsto \alpha u+ \beta v+ (z_0 + \sigma(z-z_0))/n$, where $\alpha, \beta\in \Z$ and $\sigma\in \{-1,1\}$. This tiling is built of translations of the triangle $T$, and of its symmetric with respect to the vertex $z_0$, scaled by the factor $1/n$.  %Part of the tiling $\cP_{T,n}$ appears in dark grey on Figure \ref{figAppendix}.
We define the collection of triangles 
$$
\cI_n := \{T' \in \cI_n^0 \ssep T' \subset R_n\},
$$
which is illustrated in dark gray on Figure \ref{figAppendix}. Any $T'\in \cI_n$ satisfies $\diam(T') = \diam(T)/n$ and $S(T') = S(T)$. Furthermore these elements are of the form mentioned in Point 3 of the lemma, hence are included in $\overline {R\sm R_{T,n}}$. The elements of $\cI_n$ cover the set $\{y\in R \ssep n \, d(y, \partial R) \geq \diam(T)+\diam(R)\}$, which area is larger than $|R| - C/n$. We thus obtain the left part of \iref{boundT}. Observing that $|T'| = |T|/n^2$ for each $T'\in \cT'$, we obtain that $\#(\cI_n)/n^2 \to |R|/|T|$ as $n \to \infty$.

We define a collection $\cL_n^0$ of convex polygons as follows:
$$
\cL_n^0 := \{ R_n \cap T' \sep T'\in \cI_n^0 \text{ and } \interior(T') \cap \partial R_n\neq \emptyset\},
$$
where $\interior(E)$ denotes the interior of a set $E\subset \R^2$. The set $\cL_n^0$ is illustrated in white on Figure \ref{figAppendix} (center left). The elements $T'\in \cI_n^0$ such that $\interior(T') \cap \partial R_n\neq \emptyset$ are included in the set $\{y \in \R^2 \ssep n \, d(y, \partial R) \leq \diam(R) +\diam(T)\}$ which area is smaller than $C/n$, and their individual area is $|T|/n^2$. Therefore $\#(\cL_n^0) \leq C n$. The normals to faces of the elements of $\cL_n^0$ belong to a family of at most $6$ elements: the normals to the faces of $T$, and to the faces of $R$. Therefore at most 
$6\times 5$ different angles can appear in $\cL_n^0$, and we denote the largest of these by $\alpha < \cPi$.
We denote by $\cL_n$ the collection of triangles, illustrated Figure \ref{figAppendix} (center right), obtained by triangulating each element of $\cL_n^0$, which is a convex polygon with at most six faces. The angles of the triangles partitioning a convex polygon are smaller than the angles of this polygon, hence the angles of the elements of $\cL_n$ are also bounded by $\alpha$. Furthermore $\#(\cL_n) \leq 4 \#(\cL_n^0) \leq Cn$, and $\diam(T') \leq \diam(T)/n$ for each $T'\in \cL_n$.

We denote by $E_n$ the collection of $n$ equidistributed points on each edge of $R$, described in Point 2 of Lemma \ref{lemmaTile}, and 
we denote by $E'_n$ the collection of vertices of the triangles in $\cL_n$ that fall on $\partial R'_n$. 
For each point $p\in E_n$, we draw an edge between $p$ and the point of $p'\in E'_n$ which is the closest to $p$. Note that $|p - p'| \leq d(p, R_n) + \diam(T)/n \leq (\diam(R)+ \diam(T))/n$. This produces a partition $\cB_n^0$ of $\overline {R\sm R_n}$ into triangles and convex quadrilaterals, of diameter at most 
$$
\diam(R)/n + 2( \diam(R)+ \diam(T))/n \leq 3(\diam(R)+ \diam(T))/n, %(2/3)?
$$
since the distance between two consecutive points in $E_n$ is at most $\diam(R)$. 
We denote by $\cB_n$ the collection of triangles, illustrated Figure \ref{figAppendix} (right), obtained by triangulating each polygon $K\in \cB_n^0$, of vertices $\overline K\cap (E_n \cup E'_n)$. 
We have $\#(\cB_n) = \#(E_n) + \#(E'_n) \leq 3n + 3\#(\cL_n^0)\leq C n$.
In order to conclude the proof of this lemma, we only need to show that the angles of the elements of $\cB_n$ are uniformly bounded away from $\cPi$.

We consider a triangle $T'\in \cB_n$, we denote by $L$ the length of the edge of $T'$ included in $\partial R\cup \partial R_n$, and by $H$ the height of the triangle $T'$  such that $L H = 2 |T'|$. 
It follows from \iref{distToRn} that $H \geq 2 |T| / (3n\diam(T))$.
Let $L'$ be another edge of $T'$, and let $\theta$ be the angle of $T'$ between the edges $L$ and $L'$. Then
$$
2 |T'| = L L' \sin \theta = L H,
$$
hence $\sin \theta \geq \frac H {\diam(T')} \geq \frac c C$, where $c=2 |T| / (3n\diam(T))$ and $C=3(\diam(R)+ \diam(T))$, which implies that $\arcsin(\frac c C) \leq \theta \leq \pi - \arcsin(\frac c C)$. It follows that all the angles of $T'$ are smaller than $\pi - \arcsin(\frac c C)$ which concludes the proof.

\section{Proof of Theorem \ref{thPolEq2}}
\label{appenLQ}

We consider an arbitrary but fixed $r\geq 2$ and we define $s_r := \lfloor\frac r 2\rfloor +1$. It is established in Proposition 2.1 of \cite{Mi} (equivalently Proposition 2.2.1 of \cite{thesisJM}) that for any $\pi \in \H_r$ the three following properties are equivalent  
\be
\label{vanishKE}
\left [ %\Updownarrow
\begin{array}{l}
K_r(\pi) = 0, \\
\text{There exists } \alpha, \beta\in \R \text{ and } \tilde \pi \in \H_{r-s_r} \text{ such that } \pi = (\alpha x+\beta y)^{s_r} \tilde \pi,\\
\text{There exists a sequence } (\phi_n)_{n\geq 0},\  \phi_n \in \SL_2, \text{ such that } \pi\circ \phi_n\to 0.
\end{array}
\right.
\ee
In addition the following invariance property is established in Theorem 2.6.3 of \cite{thesisJM}: let $Q$ be a polynomial on $\H_r$ such that $K_r \sim \sqrt[d] {|Q|}$ where $d=\deg Q$. Then $d r/2$ is an integer and for all $\pi \in \H_r$ and all $\phi\in \M_2$
\be
\label{invQ}
Q(\pi\circ \phi) = (\det \phi)^{\frac{d r} 2}  Q(\pi).
\ee
It follows that the polynomials $(Q_k)_{0\leq k \leq d}$, defined in \iref{defQk}, satisfy for all $\pi_1,\pi_2\in \H_r$ and all $\phi\in \M_2$
%This property clearly transfers to the polynomials $(Q_k)_{0\leq k \leq d}$ defined in Equation \iref{defQk}, in such way that 
\be
\label{invQk}
Q_k(\pi_1\circ \phi, \pi_2\circ\phi) = (\det\phi)^{\frac{d r} 2} Q_k(\pi_1, \pi_2). 
\ee
We define two functions on $\H_r \times \H_r$
\be
\label{defKstar}
\KEq(\pi_1,\pi_2) := \sqrt[2d] {\sum_{0\leq k\leq r} Q_k(\pi_1,\pi_2)^2} \stext{ and } K(\pi_1,\pi_2) := \sqrt[2 \tilde d]{\tilde Q(\pi_1^2+\pi_2^2)},
\ee
where $\tilde Q$ is an homogeneous polynomial on $\H_{2r}$ such that $K_{2r}\sim \sqrt[\tilde d]{\tilde Q}$, and $\tilde d := \deg \tilde Q$. 
We show below that $K\sim K_*$ on $\H_r \times \H_r$. 
Choosing $r=m-1$ and combining this result with \iref{eqLK} concludes the proof of Theorem \ref{thPolEq2}.

Using \iref{invQk} and remarking the invariance property $\ti Q (\pi \circ \phi) = (\det \phi)^{\ti d r} Q(\pi)$, for the same reasons as \iref{invQ}, we obtain for all $\pi_1, \pi_2\in \H_r$ and all $\phi\in \M_2$
%We deduce from Equation \iref{invQk} the invariance property
\be
\label{invPair}
%\text{for all } \pi_1, \pi_2 \in \H_r \text{ and all } \phi\in \M_2(\R) , 
%\left\{
\begin{array}{rcl}
K(\pi_1\circ \phi, \pi_2\circ\phi) &=& |\det \phi|^{\frac r 2} K(\pi_1, \pi_2),\\
\KEq(\pi_1\circ \phi, \pi_2\circ\phi) &=& |\det \phi|^{\frac r 2} \KEq(\pi_1, \pi_2).
\end{array}
%\right.
\ee

Consider $\pi_1,\pi_2\in \H_r$.
The equality $K(\pi_1, \pi_2) = 0$ is equivalent to $\pi_1^2+ \pi_2^2\in \H_{2r}$ having a linear factor of multiplicity $s_{2r} = r+1$ (according to \iref{vanishKE}), which is also equivalent to $\pi_1$ and $\pi_2$ having a common linear factor of multiplicity $s_r$. 

On the other hand the equality $\KEq(\pi_1,\pi_2) = 0$ is equivalent to $Q_k(\pi_1,\pi_2) = 0$ for all $0\leq k\leq d$. This is equivalent to $K_r(u \pi_1+ v\pi_2) = 0$ for all $u,v\in \R$ (using \iref{defQk}), which means that the polynomial $u\pi_1+ v\pi_2\in \H_r$ has a linear factor of multiplicity $s_r$ for all $u,v\in \R$ (using \iref{vanishKE}). This is equivalent to $\pi_1$ and $\pi_2$ having a common linear factor of multiplicity $s_r$. 

The following properties are therefore equivalent 
\be
\label{commonRoot}
\left[\begin{array}{l}
K(\pi_1,\pi_2) = 0,\\
\KEq(\pi_1,\pi_2) = 0,\\
\text{There exists } \alpha, \beta\in \R \text{ and } \tilde \pi_1, \tilde \pi_2 \in \H_{r-s_r} \text{ such that } \pi_1 = (\alpha x+\beta y)^{s_r} \tilde \pi_1,
\text{ and } \pi_2 = (\alpha x+\beta y)^{s_r} \tilde \pi_2.\\
\end{array}\right.
\ee
Using \iref{vanishKE} and the definition \iref{defKstar} of $K$, we find that these properties are also equivalent to 
\be
\label{commonRoot2}
\left[\begin{array}{l}
K_{2r}(\pi_1^2+\pi_2^2) = 0,\\
\text{There exists a sequence } (\phi_n)_{n\geq 0}, \ \phi_n \in \SL_2, \text{ such that } (\pi_1\circ \phi_n)^2 + (\pi_2\circ \phi_n)^2 \to 0,\\
\text{There exists a sequence } (\phi_n)_{n\geq 0}, \ \phi_n \in \SL_2, \text{ such that } \pi_1\circ \phi_n \to 0 \text{ and }\pi_2\circ \phi_n \to 0.\\
\end{array}\right.
\ee
We now define the norm $\|(\pi_1,\pi_2)\| := \sup_{|u|\leq 1} |(\pi_1(u),\pi_2(u))|$ on $\H_r\times \H_r$, and the set
$$
\cF := \{(\pi_1, \pi_2)\in \H_r\times \H_r\sep \|(\pi_1, \pi_2)\| =1 \text{ and } \|(\pi_1\circ\phi, \pi_2\circ\phi)\| \geq 1 \text{ for all } \phi\in \SL_2\}.
$$ 
The set $\cF$ is compact subset of $\H_r \times \H_r$, and $K$ as well as $K_*$ do not vanish on $\cF$ according to \iref{commonRoot} and \iref{commonRoot2}. Since these functions are continuous, there exists a constant $C_0\geq 1$ such that 
\be
\label{equivKA}
C_0^{-1} K\leq K_*\leq C_0 K \text{ on } \cF.
\ee
Let $(\pi_1, \pi_2)\in \H_r\times \H_r$. If there exists a sequence $(\phi_n)_{n\geq 0}$, $\phi_n \in \SL_2$, such that $\pi_1\circ \phi_n \to 0$ and $\pi_2 \circ \phi_n \to 0$, then $K(\pi_1, \pi_2) = K(0,0) = 0$ and $\KEq(\pi_1, \pi_2) = \KEq(0,0) = 0$ using \iref{invPair} and the continuity of $K$ and $\KEq$. Otherwise, consider a sequence $(\phi_n)_{n\geq 0}$, $\phi_n\in\SL_2$, such that 
$$
\lim_{n\to \infty}\|(\pi_1\circ \phi_n, \pi_2\circ \phi_n)\| = \inf_{\phi\in \SL_2} \|(\pi_1\circ\phi, \pi_2\circ\phi)\|.
$$
By compactness there exists a pair $(\tilde \pi_1, \tilde \pi_2)\in \H_m \times \H_m$ and a subsequence $(\phi_{n_k})_{k\geq 0}$ such that  
$$
(\pi_1\circ \phi_{n_k}, \pi_2\circ \phi_{n_k})\to (\tilde \pi_1, \tilde \pi_2).
$$ 
One easily checks that $\frac{(\ti \pi_2, \ti \pi_2)}{\|(\ti \pi_2,\ti \pi_2)\|}\in \cF$.
Using \iref{invPair} we obtain
$$
\frac{K(\pi_1,\pi_2)}{K_*(\pi_1,\pi_2)} =
\lim_{n\to \infty}  \frac{K(\pi_1 \circ\phi_n,\pi_2\circ\phi_n)}{K_*(\pi_1\circ\phi_n,\pi_2\circ\phi_n)} = \frac{K(\tilde\pi_1,\tilde\pi_2)}{K_*(\tilde \pi_1,\tilde \pi_2)}.
$$
Using  \iref{equivKA} and the homogeneity of $K$ and $K_*$, we obtain that $C_0^{-1} K\leq K_*\leq C_0 K$ on $\H_r \times \H_r$ which concludes the proof.

\subsection*{Acknowledgement}
I am extremely grateful to 
my Ph.D advisor Albert Cohen for his support
in the elaboration of this paper.

%(*pour utiliser la "vraie" biblio*)
%\cite{devore_multiscale_2009}  
%\cite{alpert_class_1993}
%\bibliographystyle{plain}
%\bibliography{../../Refs/Bibliography2.bib}

\noindent
\nl
\nl
Jean-Marie Mirebeau
\nl
UPMC Univ Paris 06, UMR 7598, Laboratoire Jacques-Louis Lions, F-75005, Paris, France
\nl
CNRS, UMR 7598, Laboratoire Jacques-Louis Lions, F-75005, Paris, France
\nl
mirebeau@ann.jussieu.fr

 \end{document}